\documentclass[12pt]{article} \usepackage{amssymb,latexsym}

\def\neweq{\setcounter{equation}{0}}

\newtheorem{theorem}{Theorem}[section]
\newtheorem{proposition}[theorem]{Proposition}
\newtheorem{corollary}[theorem]{Corollary}

\newtheorem{remark}[theorem]{Remark}
\newtheorem{lemma}[theorem]{Lemma}

\newtheorem{conjecture}[theorem]{Conjecture}

\def\proof{\smallskip\noindent {\bf Proof.}\ }
\def\endproof{\hfill$\square$\medskip}
\def\Vol{\mathrm{Vol\,}}
\def\Fix{\mathrm{Fix}}

\def\asc{\mathrm{asc}}
\def\rk{\mathrm{rk}}
\def\des{\mathrm{des}}
\def\fS{\mathfrak{S}}
\def\A{{\cal A}}
\def\C{\mathbb{C}}

\def\L{{\cal L}}

\def\str{\rule[-.3cm]{0cm}{.9cm}}
\def\poin{\mathrm{Poin}}
\def\beq{\begin{equation}}
\def\eeq{\end{equation}}
\def\beas{\begin{eqnarray*}}
\def\eeas{\end{eqnarray*}}
\def\bea{\begin{eqnarray}}
\def\eea{\end{eqnarray}}

\title{\textbf{Deformations of Coxeter Hyperplane Arrangements}}

\author{
\begin{tabular}{c}
\textsc{Alexander Postnikov} \\[.1in]
\texttt{apost@math.mit.edu}
\end{tabular}
\qquad 
\begin{tabular}{c}
\textsc{Richard P. Stanley}\\[.1in]
\texttt{rstan@math.mit.edu}
\end{tabular}\\[.3in]
Department of Mathematics\\ Massachusetts Institute of Technology\\ 
Cambridge, MA 02139}

\date{\normalsize\textbf{Version of 14 April 1997}\\ \textbf{(Not
quite complete)}} 
 
\begin{document}
\maketitle

\begin{abstract}
  We investigate several hyperplane arrangements that can be viewed as
  deformations of Coxeter arrangements. In particular, we prove a
  conjecture of Linial and Stanley that the number of regions of the
  arrangement \[ x_i - x_j = 1,\qquad 1\leq i<j\leq n, \] is equal to
  the number of alternating trees on $n+1$ vertices.  Remarkably,
  these numbers have 
  several additional combinatorial interpretations in terms of binary
  trees, partially ordered sets, and tournaments.  More generally, we
  give formulae for the number of regions and the Poincar\'e
  polynomial of certain finite subarrangements of the affine Coxeter
  arrangement of type~$A_{n-1}$.  These formulae enable us to prove a
  ``Riemann hypothesis'' on the location of zeros of the Poincar\'e
  polynomial. We also consider some generic deformations of Coxeter
  arrangements of type $A_{n-1}$.
\end{abstract}

\section{Introduction}
\label{sec:intro}
\neweq

The \emph{Coxeter arrangement} of type $A_{n-1}$ is the arrangement of
hyperplanes given by
\begin{equation}
  \label{eq:intro.1}
  x_i-x_j=0,\qquad 1\leq i < j \leq n.
\end{equation}
This arrangement has $n!$ regions.  They correspond to $n!$ different
ways of ordering the sequence $x_1,\dots,x_n$.

In the paper we extend this simple, nevertheless important, result to
the case of a general class of arrangements which can be viewed as
deformations of the arrangement~(\ref{eq:intro.1}).

One special case of such deformations is the arrangement given by
\begin{equation}
  \label{eq:intro.1.5}
  x_i-x_j=1,\qquad 1\leq i < j \leq n.
\end{equation}
We will call it the \emph{Linial arrangement.}  This arrangement was
first considered by N.~Linial and S.~Ravid.  They calculated its
number of regions and the Poincar\'e polynomial for $n\leq 9$.  On the
basis of this numerical data the second author of the present paper
made a conjecture that the number of regions of~(\ref{eq:intro.1.5})
is equal to the number of \emph{alternating trees} on $n+1$ vertices
(see~\cite{S}).  A tree $T$ on the vertices $1,2,\dots,n+1$ is
alternating if the vertices in any path in $T$ alternate, i.e., form
an up-down or down-up sequence. Equivalently, every vertex is either
less than all its neighbors or greater than all its neighbors.  These
trees first appeared in~\cite{GGP}, and in~\cite{P} a formula for the
number of such trees on $n+1$ vertices was proved.  In this paper we
provide a proof of the conjecture on the number of regions of the
Linial arrangement.

In fact, we prove a more general result for \emph{truncated affine
  arrangements}, which are certain finite subarrangements of the
affine hyperplane arrangement of type~$\widetilde{A}_{n-1}$ (see
Section~\ref{sec:affine}).  As a byproduct we get an amazing theorem
on the location of zeros of Poincar\'e polynomials of these
arrangements.  This theorem says that in one case all zeros are real,
whereas in the other case all zeros have the same real part.

The paper is organized as follows.  In Section~\ref{sec:arrangements}
we give the basic notions of hyperplane arrangement, number of
regions, Poincar\'e polynomial, and intersection poset. In
Section~\ref{sec:coxeter-arrangements} we describe the arrangements we
will be concerned with in this paper---deformations of the
arrangement~(\ref{eq:intro.1}).  In Section~\ref{sec:whitney} we
review several general theorems on hyperplane arrangements.  Then in
Section~\ref{sec:graphical} we apply these theorems to deformed
Coxeter arrangements.  In Section~\ref{sec:semigen} we consider a
``semigeneric'' deformation of the braid arrangement (the Coxeter
arrangement of type $A_{n-1}$) related to the theory of interval
orders. In Section~\ref{sec:cat} we study the hyperplane arrangements
which are related, in a special case, to interval orders
(cf.~\cite{S}) and the Catalan numbers. We prove a theorem that
establishes a relation between the numbers of regions of such
arrangements.  In Section~\ref{sec:linial} we formulate the main
result on the Linial arrangement.  We introduce several combinatorial
objects whose numbers are equal to the number of regions of the Linial
arrangement: alternating trees, local binary search trees, sleek
posets, semiacyclic tournaments.  We also prove a theorem on
characterization of sleek posets in terms of forbidden subposets.  At
last, in Section~\ref{sec:affine} we study truncated affine
arrangements.  We prove a functional equation for the generating
function for the numbers of regions of such arrangements, deduce a
formula for these numbers, and the theorem on the location of zeros of
the characteristic polynomial.

\section{Arrangements of Hyperplanes}
\neweq
\label{sec:arrangements}

First, we give several basic notions related to arrangements of
hyperplanes.  For more details, see~\cite{zas, OS, OT}.

A \emph{hyperplane arrangement} is a discrete collection of affine
hyperplanes in a vector space.  We will be concerned here only with
finite arrangements.  Let $\A$ be a finite hyperplane arrangement in a
real finite-dimensional vector space~$V$.  It will be convenient to
assume that the vectors dual to hyperplanes in $\A$ span the vector
space~$V^*$.  Denote by $r(\A)$ the number of \emph{regions} of $\A$,
which are the connected components of the space $V-\bigcup_{H\in
\A}H$.  We will also consider the number $b(\A)$ of (relatively)
\emph{bounded} regions of $\A$.

These numbers have a natural $q$-analogue.  Let $\A_{\C}$ denote the
complexified arrangement $\A$. In other words, $\A_{\C}$ is the
collection of the hyperplanes $H\otimes\C$, $H\in\A$, in the complex
vector space $V\otimes\C$.  Let $C_\A$ be the complement to
hyperplanes of~$\A_{\C}$ in $V\otimes\C$.  Then one can define the
\emph{Poincar\'e polynomial} $\poin_\A(q)$ of $\A$ as
\[
\poin_\A(q)=\sum_{k\ge 0} \dim \mathrm{H}^k(C_\A,\C)\, q^k,
\]
the generating function for the Betti numbers of $C_\A$.

The following theorem, proved in the paper of Orlik and
Solomon~\cite{OS}, shows that the Poincar\'e polynomial generalizes
the number of regions $r(\A)$ and the number of bounded regions
$b(\A)$.

\begin{theorem}
  \label{th:char} We have 
  $r(\A)=\poin_\A(1)$ and $b(\A)=\poin_\A(-1)$.
\end{theorem}

Orlik and Solomon gave a combinatorial description of the cohomology
ring $\mathrm{H}^*(C_\A,\C)$ (cf.~Section~\ref{sec:orlik-solomon}) in
terms of the \emph{intersection poset} $L_\A$ of the arrangement~$\A$.

The intersection poset is defined as follows: The elements of $L_\A$
are nonempty intersections of hyperplanes in $\A$ ordered by reverse
inclusion. The poset $L_\A$ has a unique minimal element $\hat{0}=V$.
This poset is always a meet-semilattice for which every interval is a
geometric lattice. It will be a (geometric) lattice if and only if
$L_\A$ contains a unique maximal element, i.e., the intersection of all
hyperplanes in $\A$ is nonempty.  In fact, $L_\A$ is a geometric
semilattice in the sense of Wachs and Walker~\cite{ww}, and thus for
instance is a shellable and hence Cohen-Macaulay poset.

The \emph{characteristic polynomial} of $\A$ is defined by
\begin{equation}
  \chi_{\A}(q) = \sum_{z\in L_{\A}}\mu(\hat{0},z)\,q^{\dim z},
\end{equation}
where $\mu$ denotes the M\"obius function of $L_{\A}$
(see \cite[Section~3.7]{rs:ec}).

Let $d$ be the dimension of the vector space~$V$.  Note that it follows
from the properties of geometric lattices
\cite[Proposition~3.10.1]{rs:ec} that the sign of $\mu(\hat{0},z)$ is
equal to~$(-1)^{d-\dim z}$.

The following simple relation between the (topologically defined)
Poincar\'e polynomial and the (combinatorially defined) characteristic
polynomial was found in~\cite{OS}:
\begin{equation}
  \label{eq:poinc_char}
  \chi_{\A}(q)=q^{d}\poin_\A(-q^{-1}).
\end{equation}
Sometimes it will be more convenient for us to work with the
characteristic polynomial $\chi_A(q)$ rather than the Poincar\'e
polynomial.

A combinatorial proof of Theorem~\ref{th:char} in terms of the
characteristic polynomial was earlier given by T.~Zaslavsky
in~\cite{zas}.

The number of regions, the number of (relatively) bounded regions,
and, more generally, the Poincar\'e (or characteristic) polynomial are
the most simple numerical invariants of a hyperplane arrangement.  In
this paper we will calculate these invariants for several hyperplane
arrangements related to Coxeter arrangements.

\section{Coxeter Arrangements and their Deformations}
\label{sec:coxeter-arrangements}
\neweq

Let $V_{n-1}$ denote the subspace (hyperplane) in $\mathbb{R}^n$ of
all vectors $(x_1,\dots,x_n)$ such that $x_1+\cdots+x_n=0$.  All
hyperplane arrangements that we consider below lie in $V_{n-1}$.  The
lower index $n-1$ will always denote dimension of an arrangement.

The \emph{braid arrangement} or \emph{Coxeter arrangement (of type
  $A_{n-1}$}) is the arrangement $\A_{n-1}$ of hyperplanes in
$V_{n-1}\subset\mathbb{R}^n$ given by
\begin{equation}
  \label{eq:cox}
  x_i-x_j=0,\qquad 1\le i<j\le n.
\end{equation} 


It is clear that $\A$ has $r(\A_{n-1})=n!$ regions (called Weyl
chambers) and $b(\A_{n-1})=0$ bounded regions.  Arnold~\cite{Arnold}
calculated the cohomology ring $H^*(C_{\A_n},\mathbb{C})$. In
particular, he proved that
\begin{equation}
  \label{eq:arnold}
  \poin_{\A_{n-1}}(q)=(1+q)(1+2q)\cdots (1+(n-1)q).
\end{equation}

\setlength{\unitlength}{.5pt}
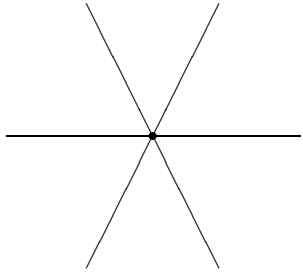
\begin{figure}
  \begin{center}
    \leavevmode
    \begin{picture}(200,200)(0,0)
      \put(50,0){\line(1,2){100}}
      \put(150,0){\line(-1,2){100}}
      \put(-11,100){\line(1,0){223}}
      \put(100,100){\circle*{6}}
    \end{picture}
    \caption{The Coxeter hyperplane arrangement $\A_2$.}
    \label{fig:A_2}
  \end{center}
\end{figure}

In this paper we will study \emph{deformations} of the
arrangement~(\ref{eq:cox}), which are hyperplane arrangements
in~$V_{n-1}\subset\mathbb{R}^n$ of the following type:

\begin{equation}
  x_i-x_j= a_{ij}^{(1)},\dots,a_{ij}^{(m_{ij})},\qquad 1\le i<j\le n.
  \label{1.1}
\end{equation}
where $m_{ij}$ are nonnegative integers and $a_{ij}^{(k)}\in
\mathbb{R}$.

One special case is the arrangement given by
\begin{equation}
  x_i-x_j=a_{ij},\qquad 1\le i< j\le n.
\label{1.05}
\end{equation}

\medskip The following hyperplane arrangements of type~(\ref{1.1})
worth mentioning:
\begin{itemize}
\item The \emph{generic arrangement} (see the end of
Section~\ref{sec:graphical}) given by
  $$ x_i - x_j = a_{ij}, \qquad 1\leq i<j\leq n, $$
where the $a_{ij}$'s are generic real numbers.
\item The \emph{semigeneric arrangement} ${\cal G}_n$ (see
Section~\ref{sec:semigen}) given by
  $$ x_i-x_j=a_i, \qquad 1\leq i\leq n,\ 1\leq j\leq n,\ i\neq j, $$
where the $a_i$'s are generic real numbers.
\item 
  The \emph{Linial arrangement} $\L_{n-1}$ (see \cite{S} and
  Section~\ref{sec:linial}) given by
  \begin{equation}
    \label{eq:1.linial}
    x_i - x_j = 1,\qquad 1\leq i<j\leq n.
  \end{equation}
\item
  The \emph{Shi arrangement} ${\cal S}_{n-1}$ (see~\cite{shi, shi2, S} 
  and Section~\ref{sec:affine-formulae}) given by
  \begin{equation}
    \label{eq:1.shi}
    x_i - x_j = 0, 1, \qquad 1\leq i<j\leq n.
  \end{equation}
\item 
  The \emph{extended Shi arrangement} ${\cal S}_{n-1,\, k}$
  (see Section~\ref{sec:affine-formulae}) 
  given by
  \begin{equation}
    \label{eq:1.ext.shi}
    x_i - x_j = -k, -k+1,\dots, k+1, \qquad 1\leq i<j\leq n,
  \end{equation}
  where $k\geq 0$ is fixed.
\item
  The \emph{Catalan arrangements} (see Section~\ref{sec:cat}) 
  ${\cal C}_{n-1}(1)$ given by
  \begin{equation}
    \label{eq:1.cat}
    x_i - x_j = -1,1,\qquad 1\leq i<j\leq n,
  \end{equation}
  and ${\cal C}_{n-1}^0(1)$ given by
  \begin{equation}
    \label{eq:1.cat.0}
    x_i - x_j = -1,0,1,\qquad 1\leq i<j\leq n.
  \end{equation}
\item
  The \emph{truncated affine arrangement} $\A_{n-1}^{ab}$ 
  (see Section~\ref{sec:affine})
  given by
  \begin{equation}
    \label{eq:1.truncate}
    x_i - x_j = -a+1, -a+2,\dots, b-1 ,\qquad 1\leq i<j\leq n,
  \end{equation}
  where $a$ and $b$ are fixed integers such that $a+b\geq 2$.
\end{itemize}

\setlength{\unitlength}{.5pt}
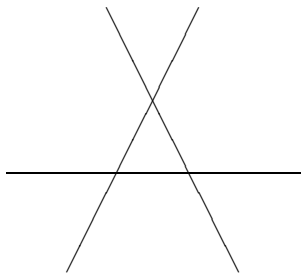
\begin{figure}
  \begin{center}
    \leavevmode
    \begin{picture}(200,200)(0,0)
      \put(35,0){\line(1,2){100}}
      \put(165,0){\line(-1,2){100}}
      \put(-11,75){\line(1,0){223}}
    \end{picture}
    \caption{Seven regions of the Linial arrangement $\L_2$.}
    \label{fig:L_3}
  \end{center}
\end{figure}

\medskip One can define analogous arrangements for any root system.
Let $V$ be a real $d$-dimensional vector space, and let $R$ be a
\emph{root system} in $V^*$ with a chosen set of \emph{positive} roots
$R_+=\{\beta_1,\beta_2,\dots,\beta_N\}$ (see,
e.g.,~\cite[Ch.~VI]{Bu}).  The \emph{Coxeter arrangement}~${\cal R}$ of
type $R$ is the arrangement of hyperplanes in $V$ given by
\begin{equation}
  \label{eq:R-coxeter}
  \beta_i(x)=0,\qquad 1\leq i \leq N.
\end{equation}

Brieskorn~\cite{brieskorn} generalized Arnold's
formula~(\ref{eq:arnold}).  His formula for the Poincar\'e polynomial
of~(\ref{eq:R-coxeter}) involves the exponents $e_1,\dots,e_d$ of the
corresponding Weyl group~$W$:
\[
\poin_{\cal R}(q)=(1+e_1q)(1+e_2q)\cdots (1+e_dq).
\]

Consider the hyperplane arrangement given by
\begin{equation}
\beta_i(x)=a_i^{(1)},\dots,a_i^{(m_i)}\qquad 1\le i\le N,
\label{1.2}
\end{equation}
where $x\in V$, $m_i$ are some nonnegative integers, and $a_i^{(k)}\in
\mathbb{R}$.  Many of the results of this paper have a natural
counterpart in the case of an arbitrary root system.  We will briefly
outline several related results and conjectures.  In more detail they
will appear elsewhere.

\section{Whitney's formula and the NBC theorem}
\neweq
\label{sec:whitney}

In this section we review several essentially well-known results on
hyperplane arrangements that will be useful in the what follows.

Consider the arrangement $\A$ of hyperplanes in $V\cong\mathbb{R}^d$
given by equations
\begin{equation}
\label{eq:arb}
h_i(x)=a_i, \qquad 1\le i\le N,
\end{equation}
where $x\in V$, the $h_i\in V^*$ are linear functionals on $V$, and
the $a_i$ are real numbers.

We call a subset $I$ in $\{1,2,\dots,N\}$ \emph{central} if the
intersection of the hyperplanes $h_i(x)=a_i$, $i\in I$, is nonempty.
For a subset $I=\{i_1,i_2,\dots,i_l\}$, denote by $\rk (I)$ the
dimension (rank) of the linear span of the vectors
$h_{i_1},\dots,h_{i_l}$.

The following statement is a generalization of a classical formula of
Whitney~\cite{whitney}.

\begin{theorem}
\label{th:alt_sum}
The Poincar\'e and characteristic polynomials of the
arrangement~$\A$ are equal to
\begin{eqnarray}
\label{eq:alt_sum}
\poin_\A(q) & = & \sum_I(-1)^{|I|-\rk(I)}\, q^{\rk(I)},\\[.05in]
\label{eq:alt_sum_char}
\chi_\A(q) & = & \sum_I (-1)^{|I|} \, q^{d-\rk(I)},
\end{eqnarray}
where $I$ ranges over all central subsets in $\{1,2,\dots,N\}$.  In
particular,

\begin{eqnarray}
r(\A)&=& \displaystyle\sum_I (-1)^{|I|-\rk(I)} \label{eq:whitr}\\[.2in] 
b(\A)&=& \displaystyle\sum_I (-1)^{|I|}. \nonumber \end{eqnarray}
\end{theorem}

We also need the well-known cross-cut theorem (see,
\cite[Corollary~3.9.4]{rs:ec}).
\begin{theorem}
  \label{th:cross-cut}
  Let $L$ be a finite lattice with minimal element~$\hat{0}$ and
  maximal element~$\hat{1}$, and let $X$ be a subset of vertices
  in $L$ such that {\rm (a)} $\hat{0}\not\in X$, and {\rm (b)} if $y\in
  L$ and $y\ne \hat{0}$, then $x\leq y$ for some $x\in X$ (such elements
  are called atoms).  Then
  \begin{equation}
  \label{eq:cross-cut}
    \mu_L(\hat{0},\hat{1})=\sum_k (-1)^k \, n_k,
  \end{equation}
  where $n_k$ is the number of $k$-element subsets in $X$ with
  join equal to $\hat{1}$.
\end{theorem}

Now we can easily deduce Theorem~\ref{th:alt_sum}.

\proof Let $z$ be any element in the intersection poset $L_\A$, and
let $L(z)$ be the subposet of all elements $x\in L_\A$ such that
$x\leq z$, i.e., the subspace~$x$ contains~$z$. In fact, $L(z)$ is a
geometric lattice. Let $X$ be the set of all hyperplanes from $\A$
which contain $z$. If we apply Theorem~\ref{th:cross-cut} to $L=L(z)$
and sum~(\ref{eq:cross-cut}) over all $z\in L_\A$, we get the
formula~(\ref{eq:alt_sum_char}).  Then, by~(\ref{eq:poinc_char}), we
get~(\ref{eq:alt_sum}).  \endproof

A \emph{cycle} is a minimal subset $I$ such that $\rk(I)=|I|-1$.
In other words, a subset $I=\{i_1,i_2,\dots,i_l\}$ is a cycle if there
exists a nonzero vector $(\lambda_1,\lambda_2,\dots,\lambda_l)$,
unique up to a nonzero factor, such that
$\lambda_1h_{i_1}+\lambda_2h_{i_2}+\cdots+\lambda_lh_{i_l}=0$.  It is
not difficult to see that a cycle $I$ is central if, in addition, we
have $\lambda_1a_{i_1}+\lambda_2a_{i_2}+\cdots+\lambda_l a_{i_l}=0$.
Thus, if $a_1=\cdots=a_N=0$ then all cycles are central, and if the
$a_i$ are generic then there are no central cycles.

A subset $I$ is called \emph{acyclic} if $|I|=\rk(I)$, i.e.,
$I$ contains no cycles. It is clear that any acyclic subset is central.

\begin{corollary}
  \label{cor:generic}
  In the case when the $a_i$ are generic,
  the Poincar\'e polynomial is given by
  \[
  \poin_\A(q) =\sum_I q^{\rk(I)},
  \]
  where the sum is over all acyclic subsets $I$ of $\{1,2,\dots,N\}$.
  In particular,  the number of regions $r(\A)$ is equal to
  the number of acyclic subsets.
\end{corollary}

Indeed, in this case a subset $I$ is acyclic if and only if it is
central.

\begin{remark}
{\rm The word ``generic'' in the corollary means that no~$k$ distinct
hyperplanes in~(\ref{eq:arb}) intersect in an affine subspace of
codimension less than $k$.  For example, if ${\cal A}$ is defined over
$\mathbb{Q}$ then it is sufficient to require that the $a_i$ be
linearly independent over $\mathbb{Q}$.}
\label{rem2.2}
\end{remark}

Let us fix a linear order $\rho$ on the set $\{1,2,\dots,N\}$.  We say
that a subset $I$ in $\{1,2,\dots,N\}$ is a \emph{broken central circuit} if
there exists $i\not\in I$ such that $I\cup \{i\}$ is a central cycle and
$i$ is the minimal element of $I\cup \{i\}$ with respect to the
order~$\rho$.

The following, essentially well-known, theorem gives us the main tool
for the calculation of Poincar\'e (or characteristic) polynomials.  We
will refer to it as the No Broken Circuit (NBC) Theorem.
\begin{theorem} 
\label{th:broken_general}
We have 
\[
\poin_\A(q)=\sum_I q^{|I|},
\]
where the sum is over all acyclic subsets $I$ of $\{1,2,\dots, N\}$
without broken central circuits.
\end{theorem}

\proof We will deduce this theorem from Theorem~\ref{th:alt_sum} using
the \emph{involution principle.} In order to do this we construct an
involution $\iota:I\to \iota(I)$ on the set of all central subsets $I$
with a broken central circuit such that for any $I$ we have
$\rk(\iota(I))=\rk(I)$ and $|\iota\cdot I|=|I|\pm 1$.

This involution is defined as follows:
Let $I$ be a central subset with a broken central circuit, and let
$s(I)$ be the set of all $i \in {1,\dots, N}$ such that $i$ is
the minimal element of a broken central circuit $J\subset I$.
Note that $s(I)$ is nonempty.
If the minimal element $s_*$ of $s(I)$ lies in $I$, then we define
$\iota(I)=I\setminus\{s_*\}$.  Otherwise, 
we define $\iota(I) = I\cup\{s_*\}$.

Note that $s(I)=s(\iota(I))$, thus $\iota$ is indeed an
involution.  It is clear now that all terms in
(\ref{eq:alt_sum}) for $I$ with a broken central circuit cancel each
other and the remaining terms yield the formula in
Theorem~\ref{th:broken_general}.
\endproof

\begin{remark}
  {\rm Note that by Theorem~\ref{th:broken_general} the number of
    subsets $I$ without broken central circuits does not depend on the
    choice of the linear order $\rho$.}
\label{rem3.2}
\end{remark}

\section{Deformations of Graphic Arrangements}
\label{sec:graphical}
\neweq

In this section we show how to apply the results of the previous
section to arrangements of type~(\ref{1.1}) and to give an
interpretation of these results in terms of (colored) graphs.

With the hyperplane $x_i-x_j=a_{ij}^{(k)}$ in~(\ref{1.1})
one can associate the edge $(i,j)$ that has the color~$k$. 
We will denote this edge by $(i,j)^{(k)}$.  Then a
subset $I$ of hyperplanes corresponds to a colored graph $G$ on the
set of vertices $\{1,2,\dots,n\}$. According to the definitions in
Section~\ref{sec:whitney}, a circuit $(i_1,i_2)^{(k_1)},
(i_2,i_3)^{(k_2)},\dots,(i_l,i_1)^{(k_l)}$ in $G$ is 
\emph{central} if $a_{i_1,i_2}^{(k_1)}+ a_{i_2,i_3}^{(k_2)}+\cdots+
a_{i_l,i_1}^{(k_l)}=0$. 
Clearly, a graph $G$ is acyclic if and only if $G$ is a \emph{forest}.

Fix a linear order on the edges $(i,j)^{(k)}$, $1\le i<j\le n$,
$1\le k\le m_{ij}$. 
We will call a subset of edges $C$ 
a \emph{broken $A$-circuit} if $C$ is obtained from a central circuit
by deleting the minimal element (here $A$ stands for the collection
$\{a_{ij}^{(k)}\}$). Note that it should not be confused with
the classical notion of a broken circuit of a graph, which 
corresponds to the case when all $a_{ij}^{(k)}$ are zero.

We summarize below several special cases of the NBC Theorem
(Theorem~\ref{th:broken_general}).
Here $|F|$ denotes the number of edges in a forest $F$.

\begin{corollary}
\label{cor:broken_graphic}
The Poincar\'e polynomial of the arrangement~(\ref{1.1}) is equal to
\[
\poin_\A(q)=\sum_F q^{|F|},
\]
where the sum is over all colored forests $F$ on the vertices
$1,2,\dots,n$ (an edge $(i,j)$ can have a color $k$, where$1\le k\le
m_{ij}$) without broken $A$-circuits.
The number of regions of
arrangement~(\ref{1.1}) is equal to the number of such forests.
\end{corollary}

In the case of the arrangement~(\ref{1.05}) we have:
\begin{corollary}
The Poincar\'e polynomial of the arrangement~(\ref{1.05}) is equal to
\[
\poin_\A(q)=\sum_F q^{|F|},
\]
where the sum is over all forests on the set of vertices
$\{1,2,\dots,n\}$ without broken $A$-circuits.
The number of regions of the arrangement~(\ref{1.05}) is equal to the number
of such forests.
\end{corollary}

In the case when the $a_{ij}^{(k)}$ are generic these 
results become especially simple.

For a forest $F$ on vertices $1,2,\dots,n$ we will write
$m^F:=\prod_{(i,j)\in F}m_{ij}$, where the product is over all edges
$(i,j)$, $i<j$, in $F$. Let $c(F)$ denote the number of connected
components in~$F$.

\begin{corollary} Fix nonnegative integers $m_{ij}$, $1\le i<j\le n$.
Let $\A$ be an arrangement of type~(\ref{1.1}) where
the $a_{ij}^{(k)}$ are generic. Then
\begin{enumerate} 
\item $\poin_\A(q)  =  \sum_F m^F q^{|F|}$,
\item $r(\A) = \sum_F m^F$,
\end{enumerate}
where the sums are over all forests $F$ on the vertices $1,2,\dots,n$.
\label{th2.1}
\end{corollary}

\begin{corollary} The number of regions of the arrangement~(\ref{1.05}) 
with generic $a_{ij}$ is equal to the number of forests on $n$
labelled vertices.
\label{cor2.3}
\end{corollary}

This corollary is ``dual'' to the following known result (see,
e.g.,~\cite[Exercise~4.32(a)]{rs:ec}).

\begin{proposition}
  Let $P_n$ be the permutohedron, i.e., the polyhedron with 
  vertices $(\sigma_1,\dots,\sigma_n)\in\mathbb{R}^n$, where 
  $\sigma_1,\dots,\sigma_n$ ranges over all permutations of
  $1,\dots,n$. 
  Then 
  the number of integer points in $P_n$ is 
  equal to the number of forests on $n$ vertices. 
\end{proposition}

The connected components of the $n\choose 2$-dimensional
space of all arrangements~(\ref{1.05}) correspond to 
(coherent) zonotopal tilings of the permutohedron $P_n$,
i.e., certain subdivisions of $P_n$ into parallelopipeds. 
The regions of a generic arrangement~(\ref{1.05}) correspond to 
the vertices of the corresponding tiling, which are all integer
points in~$P_n$.

\section{A semigeneric deformation of the braid arrangement.}
\label{sec:semigen}
\neweq

Define the ``semigeneric'' deformation ${\cal G}_n$ of the braid
arrangement (\ref{eq:cox}) to be the arrangement
  $$ x_i-x_j=a_i, \qquad 1\leq i\leq n,\ 1\leq j\leq n,\ i\neq j, $$
where the $a_i$'s are generic real numbers (e.g., linearly independent
over $\mathbb{Q}$). The significance of this arrangement to the theory
of interval orders is discussed in \cite[{\S}3]{S}. In \cite[Thm.\ 3.1
and Cor.\ 3.3]{S} a generating function for the number $r({\cal G}_n)$
of regions and for the characteristic polynomial $\chi_{{\cal
G}_n}(q)$ of ${\cal G}_n$ is stated without proof. In this section we
provide the proofs.

\begin{theorem}
  \label{thm:semigen}
Let
  \begin{eqnarray*} z & = & \sum_{n\geq 0} r({\cal
    G}_n)\frac{x^n}{n!}\\ & = & 1+x+3\frac{x^2}{2!}+19\frac{x^3}{3!}
    +195\frac{x^4}{4!}+2831\frac{x^5}{5!}+53703\frac{x^6}{6!}
    +\cdots. \end{eqnarray*}
Define a power series 
  $$ y = 1 + x +5\frac{x^2}{2!}+46\frac{x^3}{3!}+631\frac{x^4}{4!}
     +11586\frac{x^5}{5!}+\cdots $$ 
by the equation
  $$ 1 = y(2-e^{xy}). $$
Then $z$ is the unique power series satisfying
  $$ \frac{z'}{z} = y^2, \qquad z(0)=1. $$
\end{theorem}

\textbf{Proof.} We use the formula (\ref{eq:whitr}) to compute
$R({\cal G}_n)$. Given a central set $I$ of hyperplanes $x_i-x_j=a_i$
in ${\cal G}_n$, define a directed graph $G_I$ on the vertex set
$1,2,\dots,n$ as follows: let $i\rightarrow j$ be a directed edge of
$G_I$ if and only if the hyperplane $x_i-x_j=a_i$ belongs to $I$. (By
slight abuse of notation, we are using $I$ to denote a set of
hyperplanes, rather than the set of their indices.) Note that $G_I$
cannot contain both the edges $i\rightarrow j$ and $j\rightarrow i$,
since the intersection of the corresponding hyperplanes is empty. If
$k_1, k_2,\dots, k_r$ are distinct elements of $\{1,2,\dots,n\}$, then
it is easy to see that if $r$ is even then there are exactly two ways
to direct the edges $k_1k_2, k_2k_3, \dots, k_{r-1}k_r, k_rk_1$ so
that the hyperplanes corresponding to these edges have nonempty
intersection, while if $r$ is odd then there are no ways. It follows
that $G_I$, ignoring the direction of edges, is bipartite (i.e., all
circuits have even length). Moreover, given an undirected bipartite
graph on the vertices $1,2,\dots,n$ with blocks (maximal connected
subgraphs that remain connected when any vertex is removed) $B_1,
\dots, B_s$, there are exactly two ways to direct the edges of each
block so that the resulting directed graph $G$ is the graph $G_I$ of a
\emph{central} set $I$ of hyperplanes. In addition, rk$(I) = n-c(G)$,
where $c(G)$ is the number of connected components of $G$. Letting
$e(G)$ be the number of edges and $b(G)$ the number of blocks of $G$,
it follows from equation (\ref{eq:alt_sum_char}) that
  $$ \chi_{{\cal G}_n}(q) = \sum_G (-1)^{e(G)}2^{b(G)}q^{c(G)}, $$
where $G$ ranges over all bipartite graphs on the vertex set
$1,2,\dots, n$. This formula appears without proof in \cite[Thm.\
3.2]{S}. In particular, putting $q=-1$ gives
  \beq r({\cal G}_n) = (-1)^n\sum_G (-1)^{e(G)+c(G)}
    2^{b(G)}. \label{eq:rgn} \eeq

To evaluate the generating function $z=\sum r({\cal
G}_n)\frac{x^n}{n!}$, we use the following strategy.

\begin{enumerate}
  \item[(a)] Compute $A_n:=\sum_G (-1)^{e(G)}$, where $G$ ranges over
  \emph{all} (undirected) bipartite graphs on $1,2,\dots,n$.
  \item[(b)] Use (a) and the exponential formula to compute $B_n:=
  \sum_G (-1)^{e(G)}$, where now $G$ ranges over all \emph{connected}
  bipartite graphs on $1,2,\dots,n$. 
  \item[(c)] Use (b) and the block-tree theorem to compute the sum
  $C_n:= \sum_G (-1)^{e(G)}$, where $G$ ranges over all bipartite
  \emph{blocks} on $1,2,\dots,n$. 
  \item[(d)] Use (c) and the block-tree theorem to compute the sum
  $D_n:= \sum_G (-1)^{e(G)}2^{b(G)}$, where $G$ ranges over all
  \emph{connected} bipartite graphs on $1,2,\dots,n$.
  \item[(e)] Use (d) and the exponential formula to compute the
  desired sum (\ref{eq:rgn}).
\end{enumerate}

We now proceed to steps (a)--(e).
\vspace{.1in}

(a) Let $b_k(n)$ be the number of $k$-edge bipartite graphs on the
vertex set $1,2,\dots,n$. It is known (e.g., \cite[Exercise 5.5]{ec2})
that
  $$ \sum_{n\geq 0} \sum_{k\geq 0} b_k(n) q^k\frac{x^n}{n!} = 
   \left[ \sum_{n\geq 0}\left( \sum_{i=0}^n (1+q)^{i(n-i)}
    {n\choose i}\right)\frac{x^n}{n!}\right]^{1/2}. $$
Put $q=-1$ to get
  $$ \sum_{n\geq 0} A_n\frac{x^n}{n!} = \left( 1+\sum_{n\geq 1}
    2\frac{x^n}{n!}\right)^{1/2} = \left( 2e^x-1\right)^{1/2}. $$

(b) According to the exponential formula \cite[p.\ 166]{GJ}, we have
  \beas \sum_{n\geq 1} B_n\frac{x^n}{n!} & = & \log \sum_{n\geq 0}
  A_n\frac{x^n}{n!}\\ & = & \frac 12 \log(2e^x-1). \eeas

(c) Let $B'_n$ denote the number of \emph{rooted} connected bipartite
graphs on $1,2,\dots,n$. Since $B'_n=nB_n$, we get
  \bea \sum_{n\geq 1} B'_n\frac{x^n}{n!} & = & x\frac{d}{dx}
    \sum_{n\geq 1} B_n\frac{x^n}{n!}\nonumber \\ & = &
    \frac{x}{2-e^{-x}}. \label{eq:stepc} \eea
Suppose now that ${\cal B}$ is a set of nonisomorphic blocks $B$ and
$w$ is a weight function on ${\cal B}$, so $w(B)$ denotes the weight
of the block $B$. Let
  $$ T(x) = \sum_{B\in {\cal B}} w(B)\frac{x^{p(B)}}{p(B)!}, $$
where $p(B)$ denotes the number of vertices of $B$. Let
  $$ u(x) = \sum_G \left(\prod_B w(B)\right) \frac{x^{p(G)}}{p(G)!}, $$
where $G$ ranges over all connected graphs whose blocks are rooted and
are isomorphic (as unrooted graphs) to elements of ${\cal B}$, and
where $B$ ranges over all blocks of $G$. The \emph{block-tree theorem}
\cite[(1.3.3)]{h-p}\cite[Ch.\ 5 Exercises]{ec2} asserts that
  \beq u = xe^{T'(u)}. \label{eq:uxbu} \eeq
If we take ${\cal B}$ to be the set of all nonisomorphic bipartite
blocks, $w(B)=(-1)^{e(B)}$, and $u=x/(2-e^{-x})$, then it follows
from (\ref{eq:stepc}) that
  \beq T(x) = \sum_{n\geq 1} C_n\frac{x^n}{n!}. \label{eq:bt1} \eeq

(d) Let $D'_n$ be defined like $D_n$, except that $G$ ranges over all
\emph{rooted} connected bipartite graphs on $1,2,\dots, n$, so
$D'_n=nD_n$. 
Let $v(x) = \sum_{n\geq 1} D'_n\frac{x^n}{n!}$. By the block-tree
theorem we have 
  $$ v = xe^{2T'(v)}, $$
where $T(x)$ is given by (\ref{eq:bt1}). Substitute $v^{\langle
-1\rangle}$ for $x$ and use (\ref{eq:uxbu}) to get
  \beas x & = & v^{\langle -1\rangle}(x)e^{2T'(x)}\\ & = &
    v^{\langle -1\rangle}(x)\left( \frac{x}{u^{\langle -1\rangle}(x)}
    \right)^2. \eeas
Substitute $v(x)$ for $x$ to obtain
  $$ x\, v(x) = u^{\langle -1\rangle}(v(x))^2. $$
Take the square root of both sides and compose with $u(x) =
x/(2-e^{-x})$ on the left to get
  \beq \frac{\sqrt{xv}}{2-e^{-\sqrt{xv}}} = v. \label{eq:stepd} \eeq
    
(e) Equation (\ref{eq:rgn}) and the exponential formula show that 
  \bea z & = & \exp\left( -\sum_{n\geq 1} (-1)^nD_n\frac{x^n}{n!}
   \right)\nonumber \\ & = & \exp\left(-\int \frac{v(-x)}{x}\right), 
   \label{eq:z} \eea
where $\int$ denotes the formal integral, i.e., $\int\sum
a_n\frac{x^n}{n!} = \sum a_n\frac{x^{n+1}}{(n+1)!}$. (The first minus
sign in (\ref{eq:z}) corresponds to the factor $(-1)^{c(G)}$ in
(\ref{eq:rgn}).) 

Let $v(-x)=-xy^2$. Equation (\ref{eq:stepd}) becomes (taking care
to choose the right sign of the square root)
  $$ 1 = y(2-e^{xy}), $$
while (\ref{eq:z}) shows that $z'/z = -v(-x)/x = y^2$. This completes
the proof. $\ \Box$

\textsc{Note.} The semigeneric arrangement ${\cal G}_n$ satisfies the
hypotheses of \cite[Thm.\ 1.2]{S}. It follows that 
  $$ \sum_{n\geq 0} \chi_{{\cal G}_n}(q)\frac{x^n}{n!} = z(-x)^{-q},
   $$
as stated in \cite[Cor.\ 3.3]{S}. Here $z$ is as defined in
Theorem~\ref{thm:semigen}. 

An arrangement closely related to ${\cal G}_n$ is given by
  $$ {\cal G}'_n: \quad x_i-x_j=a_i, \qquad 1\leq i<j\leq n, $$
where the $a_i$'s are generic. The analogue of equation (\ref{eq:rgn})
is 
  $$ r({\cal G}'_n) = (-1)^n\sum_G (-1)^{e(G)+c(G)}
    2^{b(G)}, $$
where now $G$ ranges over all bipartite graphs on the vertex set $1,2,
\dots, n$ for which every block is \emph{alternating}, i.e., every
vertex is either less that all its neighbors or greater than all its
neighbors. We don't see, however, how to use this formula to obtain a
generating function for $r({\cal G}'_n)$ analogous to
Theorem~\ref{thm:semigen}. 

\section{Catalan Arrangements and Semiorders}
\label{sec:cat}
\neweq

Let us fix distinct real numbers $a_1,a_2,\dots,a_m > 0$, and let
$A=(a_1,\dots,a_m)$.  In this section we consider the arrangement
${\cal C}_{n-1}={\cal C}_{n-1}(A)$ of hyperplanes in the space
$V_{n-1}= \{(x_1,\dots,x_n)\in\mathbb{R}^n\mid x_1+\cdots+x_n=0\}$
given by
\begin{equation}
\label{eq:cat1}
x_i-x_j=a_1,a_2,\dots,a_m,\quad i\ne j.
\end{equation}
We consider also the arrangement ${\cal C}_{n-1}^0={\cal C}_{n-1}^0(A)$
obtained from ${\cal C}_{n-1}$ by adjoining the hyperplanes $x_i=x_j$,
i.e., ${\cal C}_n^0$ is given by
\begin{equation}
\label{eq:cat2}
x_i-x_j=0,a_1,a_2,\dots,a_m,\quad i\ne j.
\end{equation}

Let 
\[
\begin{array}{rcl}
f_A(t)&=&\displaystyle\sum_{n\ge 0}r({\cal C}_{n-1}){ t^n\over n!},\\[.3in]
g_A(t)&=&\displaystyle\sum_{n\ge 0}r({\cal C}_{n-1}^0){t^n\over n!}
\end{array}
\]
be the exponential generating functions for the numbers of regions of
the arrangements ${\cal C}_{n-1}$ and ${\cal C}_{n-1}^0$.

The main result of this section is the following:
\begin{theorem}
\label{th:cat1}
We have $f_A(t)=g_A(1-e^{-t})$ or, equivalently,
\[
r({\cal C}_{n-1}^0)=\sum_{k\ge 0}c(n,k)\,r({\cal C}_{k-1}),
\]
where $c(n,k)$ is the signless Stirling number of the first kind,
i.e., the number of permutations of $1,2,\dots,n$ with k cycles.
\end{theorem}

Let us have a closer look at two special cases of
arrangements~(\ref{eq:cat1}) and~(\ref{eq:cat2}).  Consider the
arrangement of hyperplanes in $V_{n-1}\subset\mathbb{R}^n$ given by
the equations
\begin{equation}
\label{eq:cat3}
x_i-x_j=\pm 1,\qquad 1\leq i<j\leq n.
\end{equation}
Consider also the arrangement given by
\begin{equation}
\label{eq:cat4}
x_i-x_j=0,\ \pm 1,\qquad 1\leq i<j\leq n.
\end{equation}

It is not difficult to check the following result directly from
the definition.
\begin{proposition}
  \label{prop:cat2}
  The number of regions of the arrangement~(\ref{eq:cat4}) is equal to
  $n!\,C_n$, where $C_n$ is the Catalan number $C_n={1\over
    n+1}{2n\choose n}$.
\end{proposition}
Theorem~\ref{th:cat1} then gives a formula for the number of regions
of the arrangement~(\ref{eq:cat3}).

Let $R$ be a region of the arrangement~(\ref{eq:cat3}), and let
$(x_1,\dots,x_n)\in R$ be any point in the region $R$.  Consider the
poset $P$ on the vertices $1,\dots,n$ such that $i>_P j$ if and only
if $x_i-x_j>1$. Clearly, distinct regions correspond to distinct
posets. The posets that can be obtained in such a way are called
\emph{semiorders}. See~\cite{S} for more results on the relation
between hyperplane arrangements and \emph{interval orders} (which are
a generalization of semiorders).

The symmetric group $\fS_n$ naturally acts on the space $V_{n-1}$ by
permuting the coordinates $x_i$.  Thus it also permutes the regions of
the arrangement~(\ref{eq:cat4}). The region $x_1<x_2<\cdots<x_n$ is
called the \emph{dominant chamber}.
Every $\fS_n$-orbit consists of $n!$
regions and has a unique representative in the dominant
chamber. It is also clear that the
regions of~(\ref{eq:cat4}) in the dominant chamber correspond to
\emph{unlabelled} (i.e., nonisomorphic) semiorders on $n$ vertices.
Hence, Proposition~\ref{prop:cat2} is equivalent to a well-known
result of Wine and Freund~\cite{wi-fr} that the number of
nonisomorphic semiorders on $n$ vertices is equal to the Catalan
number.  In the special case of the arrangements~(\ref{eq:cat3})
and~(\ref{eq:cat4}), i.e., $A=(1)$, Theorem~\ref{th:cat1} gives a
formula for the number of labelled semiorders on $n$ vertices which was
first proved by Chandon, Lemaire, and Pouget~\cite{c-l-p}.

The following theorem, due to Scott and Suppes~\cite{s-s},
presents a simple characterization of semiorders
(cf.~Theorem~\ref{th:obstr}).

\begin{theorem}
\label{th:cat3}
A poset $P$ is a semiorder if and only if it contains no induced
subposet of either of the two types shown on Figure~\ref{fig:semi}.
\end{theorem}

\setlength{\unitlength}{.6pt}
\begin{figure}
\begin{picture}(300,130)(-200,-20)
\put(0,25){\circle*{5}} \put(0,75){\circle*{5}}
\put(35,25){\circle*{5}} \put(35,75){\circle*{5}}
\put(220,0){\circle*{5}} \put(220,50){\circle*{5}}
\put(220,100){\circle*{5}} \put(255,50){\circle*{5}}
\put(0,25){\line(0,1){50}} \put(35,25){\line(0,1){50}}
\put(220,0){\line(0,1){50}} \put(220,50){\line(0,1){50}}
\end{picture}
\caption{Forbidden subposets for semiorders.}
\label{fig:semi}
\end{figure}
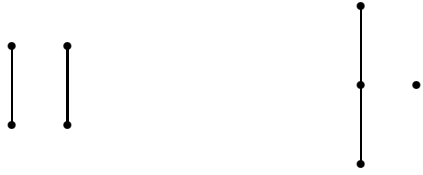

Return now to the general case of the arrangements ${\cal C}_{n-1}$ and
${\cal C}_{n-1}^0$ given by~(\ref{eq:cat1}) and~(\ref{eq:cat2}). The
symmetric group $\fS_n$ acts on the regions of ${\cal C}_{n-1}$ and ${\cal
C}_{n-1}^0$. Let $R_{n-1}$ denotes the set of all regions of ${\cal C}_{n-1}$.

\begin{lemma}
  \label{lm:cat1}
  The number of regions of ${\cal C}_{n-1}^0$ is equal to $n!$ times the
  number of $\fS_n$-orbits in $R_{n-1}$.
\end{lemma}
Indeed, the number of regions of ${\cal C}_{n-1}^0$ 
is $n!$ times the number of those in the dominant chamber.
They, in turn, correspond to $\fS_n$-orbits in $R_{n-1}$.
As was shown in~\cite{S}, the regions of ${\cal C}_{n-1}$ can be
viewed as (labelled) generalized interval orders.  On the other hand,
the regions of ${\cal C}_{n-1}^0$ that lie in the dominant chamber
correspond to unlabelled generalized interval orders.  The statement
now is tautological, that the number of unlabelled objects
is the number of $\fS_n$-orbits.

Now we can apply the following well-known lemma of Burnside (actually
first proved by Cauchy and Frobenius). 
\begin{lemma}
  \label{lm:cat2}
  Let $G$ be a finite group which acts on a finite set $M$.  Then the
  number of $G$-orbits in $M$ is equal to
  \[
  {1\over |G|}\sum_{g\in G}\Fix(g, M),
  \]
where $\Fix(g,M)$ is the number of elements in $M$ fixed by $g\in G$.
\end{lemma}

By Lemmas~\ref{lm:cat1} and~\ref{lm:cat2} we have 
\[
r({\cal C}_{n-1}^0)=\sum_{\sigma \in \fS_n}\Fix(\sigma,{\cal C}_{n-1}),
\]
where $\Fix(\sigma,{\cal C}_{n-1})$ is the number of regions of ${\cal
C}_{n-1}$ fixed by the permutation $\sigma$.

Theorem~\ref{th:cat1} now follows easily from the following lemma.
\begin{lemma}
  \label{lm:cat3}
  Let $\sigma\in \fS_n$ be a permutation with $k$ cycles.  Then the
  number of regions of ${\cal C}_{n-1}$ fixed by $\sigma$ is equal to
  the total number of regions of ${\cal C}_{k-1}$.
\end{lemma}

Indeed, by Lemma~\ref{lm:cat3}, we have
\[
r({\cal C}_{n-1}^0)=\sum_{\sigma\in \fS_n}\Fix(\sigma,{\cal C}_{n-1})=
\sum_{k\ge0}c(n,k)\,r({\cal C}_{k-1}),
\]
which is precisely the claim of Theorem~\ref{th:cat1}.

\medskip
\noindent {\bf Proof of Lemma~\ref{lm:cat3} \ }
We will construct a bijection between the regions of ${\cal C}_{n-1}$
fixed by $\sigma$ and the regions of ${\cal C}_{k-1}$.

Let $R$ be any region of ${\cal C}_{n-1}$ fixed by a permutation
$\sigma\in \fS_n$, and let $(x_1,\dots,x_n)$ be any point in $R$.  Then
for any $i,j\in\{1,\dots,n\}$ and any $s=1,\dots,m$ we have
$x_i-x_j>a_s$ if and only if $x_{\sigma(i)}-x_{\sigma(j)}>a_s$.

Let $\sigma=(c_{11}\,c_{12}\cdots c_{1l_1})\,(c_{21}\,c_{22}\cdots
c_{2l_2})\cdots (c_{k1}\,c_{k2}\cdots c_{kl_k})$ be the cycle
decomposition of the permutation
$X_i=(x_{c_{i1}},x_{c_{i2}},\dots)$ for $i=1,\dots,k$.  We will write
$X_i-X_j>a$ if $x_{i'}-x_{j'}>a$ for any $x_{i'}\in X_i$ and
$x_{j'}\in X_j$. The notation $X_i-X_j<a$ has an analogous meaning. We
will show that for any two classes $X_i$ and $X_j$ and for any
$s=1,\dots,m$ we have either $X_i-X_j>a_s$ or $X_i-X_j<a_s$.

Let $x_{i^*}$ be the maximal element in $X_i$ and let $x_{j^*}$ be the
maximal element in $X_j$.  Suppose that $x_{i^*}-x_{j^*}>a_s$. Since
$R$ is $\sigma$-invariant, for any integer $p$ we have the inequality
$x_{\sigma^p(i^*)}-x_{\sigma^p(j^*)}>a_s$.  Then, since $x_{i^*}$ is
the maximal element of $X_i$, we have $x_{i^*}-x_{\sigma^p(j^*)}>a_s$.
Again, for any integer $q$, we have
$x_{\sigma^q(i^*)}-x_{\sigma^{p+q}(j^*)}>a_s$, which implies that
$X_i-X_j>a_s$.

Analogously, suppose that $x_{i^*}-x_{j^*}<a_s$. Then for any integer
$p$ we have $x_{\sigma^p(i^*)}-x_{\sigma^p(j^*)}<a_s$.  Since
$x_{j^*}\ge x_{\sigma^p(j^*)}$, we have
$x_{\sigma^p(i^*)}-x_{j^*}<a_s$.  Finally, for any integer $q$ we
obtain $x_{\sigma^{p+q}(i^*)}-x_{\sigma^q(j^*)}<a_s$, which implies
that $X_i-X_j<a_s$.

If we pick an element $x_{i'}$ in each class $X_i$ we get a point
$(x_{1'},x_{2'},\dots,x_{k'})$ in $\mathbb{R}^k$.  This point lies in
some region $R'$ of ${\cal C}_{k-1}$. The construction above shows
that the region $R'$ does not depend on the choice of $x_{i'}$ in
$X_i$.

Thus we get a map $\phi:R\to R'$ from the regions of ${\cal C}_{n-1}$
invariant under $\sigma$ to the regions of ${\cal C}_{k-1}$.
It is clear that $\phi$ is injective. To show that $\phi$ is
surjective, let $(x_{1'},\dots,x_{k'})$ be any point in a region $R'$
of ${\cal C}_k$. Pick the point $(x_1,x_2,\dots,x_n)\in\mathbb{R}^n$
such that $x_{c_{11}}=x_{c_{12}}=\cdots=x_{1'}$,
$x_{c_{21}}=x_{c_{22}}=\cdots=x_{2'},\dots,
x_{c_{k1}}=x_{c_{k2}}=\cdots= x_{k'}$. Then $(x_1,\dots,x_n)$ is in
some region $R$ of ${\cal C}_{n-1}$ (here we use the condition
$a_1,\dots,a_m\ne 0$). According to our construction, we have
$\phi(R)=R'$.  Thus $\phi$ is a bijection.

This completes the proof of Lemma~\ref{lm:cat3} and therefore also of
Theorem~\ref{th:cat1}.  \endproof

\section{The Linial Arrangement.}
\label{sec:linial}
\neweq

As before, $V_{n-1}=\{(x_1,\dots,x_n)\in\mathbb{R}^n\mid
x_1+\cdots+x_n=0\}$.  Consider the arrangement~$\L_{n-1}$ of
hyperplanes in $V_{n-1}$ given by the equations
\begin{equation}
x_i-x_j=1,\quad 1\le i<j\le n.
\end{equation}

Recall that $r(\L_{n-1})$ denotes the number of regions of the
arrangement $\L_{n-1}$.  This arrangement was first considered by Nati
Linial and Shmulik Ravid.  They calculated the numbers $r(\L_{n-1})$
and the Poincar\'e polynomials $\poin_{\L_{n-1}}(q)$ for $n\le 9$.

In this section we give an explicit formula and several 
different combinatorial interpretations for the numbers $r(\L_{n-1})$.

\subsection{Alternating trees and local binary search trees}

We call a tree $T$ on the vertices $0,1,2,\dots,n$ \emph{alternating}
if the vertices in any path $i_1,\dots,i_k$ in $T$ alternate, i.e., we
have $i_1<i_2>i_3<\cdots\, i_k$ or $i_1>i_2<i_3>\cdots\, i_k$. In
other words, there are no $i<j<k$ such that both $(i,j)$ and $(j,k)$
are edges in~$T$. Equivalently, every vertex is either greater than
all its neighbors of less than all its neighbors.  Alternating trees
first appear in~\cite{GGP} and were studied in~\cite{P}, where they
were called \emph{intransitive trees} (see also~\cite{S}).

\setlength{\unitlength}{.8pt}
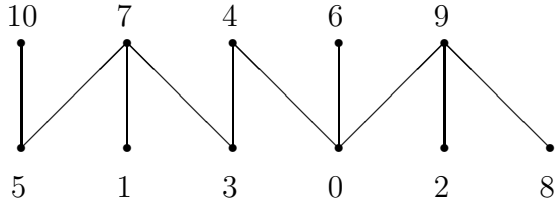
\begin{figure}[htbp]
  \begin{center}
    \leavevmode
    \begin{picture}(250,110)(0,-30)
      
      \put(0,0){\circle*{4}} 
      \put(50,0){\circle*{4}} 
      \put(100,0){\circle*{4}} 
      \put(150,0){\circle*{4}} 
      \put(200,0){\circle*{4}} 
      \put(0,50){\circle*{4}} 
      \put(50,50){\circle*{4}} 
      \put(100,50){\circle*{4}} 
      \put(150,50){\circle*{4}} 
      \put(200,50){\circle*{4}} 
      \put(250,0){\circle*{4}} 
      
      \put(0,0){\line(0,1){50}}
      \put(0,0){\line(1,1){50}}
      \put(50,0){\line(0,1){50}}
      \put(100,0){\line(-1,1){50}}
      \put(100,0){\line(0,1){50}}
      \put(150,0){\line(-1,1){50}}
      \put(150,0){\line(0,1){50}}
      \put(150,0){\line(1,1){50}}
      \put(200,0){\line(0,1){50}}
      \put(250,0){\line(-1,1){50}}

      \put(145,-23){$0$}
      \put(45,-23){$1$}
      \put(195,-23){$2$}
      \put(95,-23){$3$}
      \put(95,58){$4$}
      \put(-5,-23){$5$}
      \put(145,58){$6$}
      \put(45,58){$7$}
      \put(245,-23){$8$}
      \put(195,58){$9$}
      \put(-7,58){$10$}

    \end{picture}
    \caption{An alternating tree.}
    \label{fig:alt-tree}
  \end{center}
\end{figure}

Let $f_n$ be the number of alternating trees on the vertices 
$0,1,2,\dots,n$, and let
\[
f(x)=\sum_{n\ge 0}f_n{x^n\over n!}
\]
be the exponential generating function for the sequence $f_n$.

A plane binary tree $B$ on the vertices $1,2,\dots,n$ is called a
\emph{local binary search tree} if for any vertex $i$ in $T$ the left
child of $i$ is less than $i$ and the right child of $i$ is greater
than $i$. These trees were first considered by Ira Gessel~\cite{Ges}.
Let $g_n$ denote the number of local binary search trees on the vertices
$1,2,\dots,n$. By convention, $g_0=1$.

\setlength{\unitlength}{.6pt}
\begin{figure}[htbp]
  \begin{center}
    \leavevmode

    \begin{picture}(300,180)(0,-30)

      \put(0,0){\circle*{5}} 
      \put(100,0){\circle*{5}} 
      \put(200,0){\circle*{5}} 
      \put(300,0){\circle*{5}} 
      \put(50,50){\circle*{5}} 
      \put(150,50){\circle*{5}} 
      \put(250,50){\circle*{5}} 
      \put(100,100){\circle*{5}} 
      \put(200,100){\circle*{5}} 
      \put(150,150){\circle*{5}} 
      
      \put(0,0){\line(1,1){50}}
      \put(100,0){\line(1,1){50}}
      \put(200,0){\line(1,1){50}}
      \put(300,0){\line(-1,1){50}}
      \put(50,50){\line(1,1){50}}
      \put(150,50){\line(-1,1){50}}
      \put(250,50){\line(-1,1){50}}
      \put(100,100){\line(1,1){50}}
      \put(200,100){\line(-1,1){50}}
      
      \put(95,-23){$1$}
      \put(-5,-23){$2$}
      \put(38,55){$3$}
      \put(195,-23){$4$}
      \put(88,105){$5$}
      \put(138,155){$6$}
      \put(202,105){$7$}
      \put(152,55){$8$}
      \put(252,55){$9$}
      \put(302,5){$10$}

    \end{picture}

    \caption{A local binary search tree.}
    \label{fig:LBST}
  \end{center}
\end{figure}
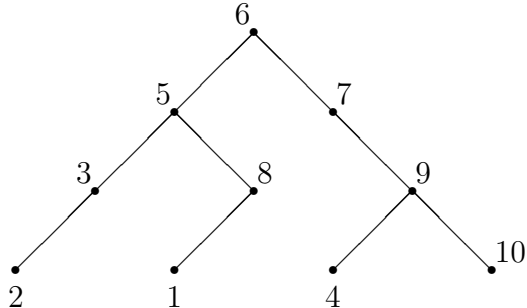

The following result was proved in~\cite{P} (see also~\cite{GGP, S}).

\begin{theorem}
For $n\ge 1$ we have 
\[
f_n=g_n=2^{-n}\sum_{k=0}^n{n\choose k}(k+1)^{n-1}
\]
and $f=f(x)$ satisfies the functional equation
\[
f=e^{x(1+f)/2}.
\]
\end{theorem}

The first few  numbers $f_n$ are given in the table below.
\medskip
\begin{center}
\begin{tabular}{|l|l|l|l|l|l|l|l|l|l|l|l|}
\hline
$n$  &\str 0 & 1 & 2 & 3 & 4  & 5   & 6    & 7     & 8      & 9       &
10\\ \hline
$f_n$&\str 1 & 1 & 2 & 7 & 36 & 246 & 2104 & 21652 & 260720 & 3598120 & 
56010096\\
\hline
\end{tabular}
\end{center}
\bigskip

      
      


The main result on the Linial arrangement is the following:

\begin{theorem}
\label{th:linial}
The number $r(\L_{n-1})$ of regions of $\L_{n-1}$ is equal to
the number $f_n$ of alternating trees on the vertices
$0,1,2\dots,n$, and thus to the number $g_n$ of local binary
search trees on $1,2,\dots,n$.
\end{theorem}

This theorem was conjectured by the second author (thanks to the
numerical data provided by Linial and Ravid) and was proved by the
first author.  A different proof was later given by C.\ 
Athanasiadis~\cite{Ath}.

In Section~\ref{sec:affine} we will prove a more general result
(see Theorems~\ref{th:af.main} and Corollary~\ref{cor:aff-ways}).

\subsection{Sleek posets and semiacyclic tournaments}
\label{sec:sleek}

Let $R$ be a region of the arrangement $\L_{n-1}$, and let
$(x_1,\dots,x_n)$ be any point in~$R$. Define $P=P(R)$ to be the poset
on the vertices $1,2,\dots,n$ such that $i<_P j$ if and only if
$x_i-x_j>1$ and $i<j$ in the usual order on $\mathbb{Z}$.

We will call a poset $P$ on the vertices $1,2,\dots,n$
\emph{sleek} if $P$ is the intersection of a semiorder
(see Section~\ref{sec:cat}) with the chain $1<2<\cdots<n$.

The following proposition immediately follows from the definitions.

\begin{proposition} The map $R\mapsto P(R)$ is 
a bijection between regions of $\L_{n-1}$ and sleek posets on
$1,2,\dots,n$.  Hence the number $r(\L_{n-1})$ is equal to the number of
sleek posets on $1,2,\dots, n$.
\end{proposition}

There is a simple characterization of sleek posets in terms of
forbidden induced subposets (compare Theorem~\ref{th:cat3}).

\begin{theorem}
\label{th:obstr}
A poset $P$ on the vertices $1,2,\dots,n$ is sleek if and only if it
contains no induced subposet of the four types shown on
Figure~\ref{fig1}, where $a<b<c<d$.
\end{theorem}

\setlength{\unitlength}{.6pt}
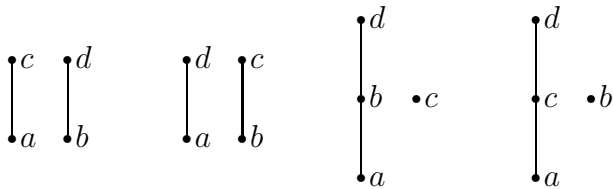
\begin{figure}\hfil
\begin{picture}(360,130)(0,-20)  
\put(0,25){\circle*{5}} \put(0,75){\circle*{5}}
\put(35,25){\circle*{5}} \put(35,75){\circle*{5}}
\put(110,25){\circle*{5}} \put(110,75){\circle*{5}}
\put(145,25){\circle*{5}} \put(145,75){\circle*{5}}
\put(220,0){\circle*{5}} \put(220,50){\circle*{5}}
\put(220,100){\circle*{5}} \put(255,50){\circle*{5}}
\put(330,0){\circle*{5}} \put(330,50){\circle*{5}}
\put(330,100){\circle*{5}} \put(365,50){\circle*{5}}
\put(0,25){\line(0,1){50}} \put(35,25){\line(0,1){50}}
\put(110,25){\line(0,1){50}} \put(145,25){\line(0,1){50}}
\put(220,0){\line(0,1){50}} \put(220,50){\line(0,1){50}}
\put(330,0){\line(0,1){50}} \put(330,50){\line(0,1){50}}
\put(5,20){$a$} \put(5,70){$c$} \put(40,20){$b$} \put(40,70){$d$}
\put(115,20){$a$} \put(115,70){$d$} \put(150,20){$b$}
\put(150,70){$c$} \put(225,-5){$a$} \put(225,45){$b$}
\put(225,95){$d$} \put(260,45){$c$} \put(335,-5){$a$}
\put(335,45){$c$} \put(335,95){$d$} \put(370,45){$b$}
\end{picture}\hfil
\caption{Obstructions to sleekness.}
\label{fig1}
\end{figure}

In the remaining part of this section we prove Theorem~\ref{th:obstr}.
 
First, we give another description of regions in $\L_{n-1}$ (or,
equivalently, sleek posets).  A \emph{tournament} on the vertices
$1,2,\dots,n$ is a directed graph $T$ without loops such that for
every $i\neq j$ either $(i,j)\in T$ or $(j,i)\in T$.  For a region $R$
of $\L_{n-1}$ construct a tournament $T=T(R)$ on the vertices
$1,2,\dots,n$ as follows: let
$(x_1,\dots,x_n)\in R$. If $x_i-x_j>1$ and $i<j$, then $(i,j)\in
T$; while if $x_i-x_j<1$ and $i<j$, then $(j,i)\in T$.

Let $C$ be a directed cycle in the complete graph $K_n$ on the
vertices $1,2,\dots,n$. We will write $C=(c_1,c_2,\dots,c_m)$ if $C$
has the edges $(c_1,c_2),(c_2,c_3),\dots,(c_m,c_1)$.  By convention,
$c_0=c_m$.  An \emph{ascent} in $C$ is a number $1\le i\le m$ such
that $c_{i-1}<c_i$.  Analogously, a \emph{descent} in $C$ is a number
$1\le i\le m$ such that $c_{i-1}>c_i$.  Let $\asc(C)$ denote the
number of ascents and $\des(C)$ denote the number of descents in $C$.
We say that a cycle $C$ is \emph{ascending} if $\asc(C)\ge \des(C)$.
For example, the following cycles are ascending: $C_0=(a,b,c)$,
$C_1=(a,c,b,d)$, $C_2=(a,d,b,c)$, $C_3=(a,b,d,c)$, $C_4=(a,c,d,b)$,
where $a<b<c<d$.  These cycles are shown on Figure~\ref{fig:asc}.

We call a tournament $T$ on $1,2,\dots, n$
\emph{semiacyclic} if it contains no ascending cycles.
In other words, $T$ is semiacyclic if for any directed cycle $C$ in
$T$ we have $\asc(C)<\des(C)$.

\setlength{\unitlength}{1.5pt}
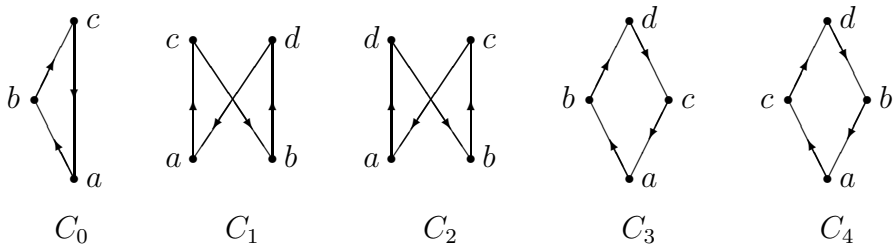
\begin{figure}\hfil
\begin{picture}(220,80)(0,-20)   

\put(0,20){\circle*{2}}
\put(10,0){\circle*{2}}
\put(10,40){\circle*{2}}
\put(5,-15){$C_0$}
\put(13,-2){$a$}
\put(13,38){$c$}
\put(-7,18){$b$}

\put(40,5){\circle*{2}}
\put(40,35){\circle*{2}}
\put(60,5){\circle*{2}}
\put(60,35){\circle*{2}}
\put(48,-15){$C_1$}
\put(33,3){$a$}
\put(63,3){$b$}
\put(33,33){$c$}
\put(63,33){$d$}

\put(90,5){\circle*{2}}
\put(90,35){\circle*{2}}
\put(110,35){\circle*{2}}
\put(110,5){\circle*{2}}
\put(98,-15){$C_2$}
\put(83,3){$a$}
\put(113,3){$b$}
\put(113,33){$c$}
\put(83,33){$d$}

\put(140,20){\circle*{2}}
\put(150,40){\circle*{2}}
\put(150,0){\circle*{2}}
\put(160,20){\circle*{2}}
\put(148,-15){$C_3$}
\put(153,-2){$a$}
\put(133,18){$b$}
\put(163,18){$c$}
\put(153,38){$d$}

\put(190,20){\circle*{2}}
\put(200,40){\circle*{2}}
\put(200,0){\circle*{2}}
\put(210,20){\circle*{2}}
\put(198,-15){$C_4$}
\put(203,-2){$a$}
\put(183,18){$c$}
\put(213,18){$b$}
\put(203,38){$d$}

\put(10,0){\line(-1,2){10}}
\put(0,20){\line(1,2){10}}
\put(10,40){\line(0,-1){40}}
\put(10,0){\vector(-1,2){5}}
\put(0,20){\vector(1,2){5}}
\put(10,40){\vector(0,-1){20}}

\put(40,5){\line(0,1){30}}
\put(40,35){\line(2,-3){20}}
\put(60,5){\line(0,1){30}}
\put(60,35){\line(-2,-3){20}}
\put(40,5){\vector(0,1){15}}
\put(40,35){\vector(2,-3){15}}
\put(60,5){\vector(0,1){15}}
\put(60,35){\vector(-2,-3){15}}

\put(90,5){\line(0,1){30}}
\put(90,35){\line(2,-3){20}}
\put(110,5){\line(0,1){30}}
\put(110,35){\line(-2,-3){20}}
\put(90,5){\vector(0,1){15}}
\put(90,35){\vector(2,-3){15}}
\put(110,5){\vector(0,1){15}}
\put(110,35){\vector(-2,-3){15}}

\put(140,20){\line(1,2){10}}
\put(150,40){\line(1,-2){10}}
\put(160,20){\line(-1,-2){10}}
\put(150,0){\line(-1,2){10}}
\put(140,20){\vector(1,2){5}}
\put(150,40){\vector(1,-2){5}}
\put(160,20){\vector(-1,-2){5}}
\put(150,0){\vector(-1,2){5}}

\put(190,20){\line(1,2){10}}
\put(200,40){\line(1,-2){10}}
\put(210,20){\line(-1,-2){10}}
\put(200,0){\line(-1,2){10}}
\put(190,20){\vector(1,2){5}}
\put(200,40){\vector(1,-2){5}}
\put(210,20){\vector(-1,-2){5}}
\put(200,0){\vector(-1,2){5}}

\end{picture}\hfil
\caption{Ascending cycles.}
\label{fig:asc}
\end{figure}

\begin{proposition}
  \label{prop:semiacyclic}
  A tournament $T$ on $1,2,\dots, n$ corresponds to a region $R$ in
  $\L_{n-1}$, i.e., $T=T(R)$, if and only if $T$ is semiacyclic.
  Hence $r(\L_{n-1})$ is the number of semiacyclic tournaments on
  $1,2,\dots, n$.
\end{proposition}
This fact was independently found by Shmulik Ravid.

For any tournament $T$ on $1,2,\dots,n$ without cycles of type
$C_0$ we can construct a poset $P=P(T)$ such that $i<_P j$ if and only
if $i<j$ and $(i,j)\in T$.  Now the four ascending cycles $C_1$, $C_2$,
$C_3$, $C_4$ in Figure~\ref{fig:asc} correspond to the four posets on
Figure~\ref{fig1}.  Therefore, Theorem~\ref{th:obstr} is equivalent to
the following result.

\begin{theorem}
  \label{th:semiacyclic}
  A tournament $T$ on the vertices $1,2,\dots,n$ is semiacyclic if and
  only if it contains no ascending cycles of the types $C_0,\ C_1,\ 
  C_2,\ C_3$, and $C_4$ shown in Figure~\ref{fig:asc}, where
  $a<b<c<d$.
\end{theorem}

\begin{remark}
{\rm This theorem is an analogue of a well-known fact that
a tournament $T$ is acyclic if and only if it contains no cycles of
length $3$. For semiacyclicity we have obstructions of lengths $3$ and
$4$.}
\end{remark} 

\proof Let $T$ be a tournament on $1,2,\dots,n$.  Suppose that $T$ is
not semiacyclic. We will show that $T$ contains a cycle of type $C_0,\ 
C_1,\ C_2,\ C_3$, or $C_4$.  Let $C=(c_1,c_2,\dots,c_m)$ be an
ascending cycle in $T$ of minimal length.  If $m=3$, or~$4$ then $C$
is of type $C_0,\ C_1,\ C_2,\ C_3$, or $C_4$.  Suppose that $m>4$.

\begin{lemma}
We have $\asc(C)=\des(C)$.

\end{lemma}

\proof
Since $C$ is ascending, we have $\asc(C)\ge \des(C)$.  Suppose
$\asc(C)>\des(c)$. If $C$ has two adjacent ascents $i$ and $i+1$ then
$(c_{i-1}, c_{i+1})\in T$ (otherwise we have an ascending cycle
$(c_{i-1},c_i,c_{i+1})$ of type $C_0$ in $T$).  Then
$C'=(c_1,c_2,\dots,c_{i-1},c_{i+1},\dots,c_m)$ is an ascending cycle
in $T$ of length $m-1$, which contradicts the fact that we chose
$C$ to be minimal.  So for every ascent $i$ in $C$ the index $i+1$ is
a descent.  Hence $\asc(C)\le \des(C)$, and we get a contradiction.
\endproof

We say that $c_i$ and $c_j$ are \emph{on the same level} in $C$ if
the number of ascents between $c_i$ and~$c_j$ is equal to the number
of descents between $c_i$ and~$c_j$.

\begin{lemma}
\label{lem:2}
We can find $i, j\in \{1,2,\dots,m\}$ such that {\rm (a)} $i$ is an
ascent and $j$ is a descent in $C$, {\rm (b)} $i\not\equiv j\pm 1\!\!
\pmod m$, and {\rm (c)} $c_{i}$ and $c_{j-1}$ are on the same level
(see Figure~\ref{fig:cyc}).
\end{lemma}

\proof
We may assume that for any $1\le s\le m$ the number of ascents in
$\{1,2,\dots,s\}$ is greater than or equal to the number of descents
in $\{1,2,\dots,s\}$ (otherwise take some cyclic permutation of
$(c_1,c_2,\dots,c_m)$).  Consider two cases.
\smallskip

\noindent 1.\quad
There exists $1\le t\le m-1$ such that $c_t$ and $c_m$ are on the same
level.  In this case, if the pair $(i,j)=(1,t)$ does not satisfy
conditions (a)--(c) then $t=2$. On the other hand, if the pair
$(i,j)=(t+1,m)$ does not satisfy (a)--(c) then $t=m-2$. Hence, $m=4$
and $C$ is of type $C_1$ or $C_2$ shown in Figure~\ref{fig:asc}.
\smallskip

\noindent 2.\quad
There is no $1\le t\le m-1$ such that $c_t$ and $c_m$ are on the same
level.  Then $2$ is an ascent and $m-1$ is a descent.  If the pair
$(i,j)=(2,m-2)$ does not satisfy (a)--(c) then
$m=4$ and $C$ is of type $C_3$ or $C_4$ shown on Figure~\ref{fig:asc}.
\endproof

Now we can complete the proof of Theorem~\ref{th:semiacyclic}.  Let
$i,j$ be two numbers satisfying the conditions of Lemma~\ref{lem:2}.
Then $c_{i-1}$, $c_i$, $c_{j-1}$, $c_j$ are four distinct vertices
such that (a) $c_{i-1}<c_i$, (b) $c_{j-1}>c_j$, (c) $c_i$ and
$c_{j-1}$ are on the same level, and (d) $c_{i-1}$ and $c_j$ are on the same
level (see Figure~\ref{fig:cyc}).  We may assume that $i<j$.

\setlength{\unitlength}{1.5pt}
\begin{figure}[htbp]
  \begin{center}
    \leavevmode
    \begin{picture}(40,60)(3,0)
      \put(-10,20){\circle*{2}}
      \put(-10,40){\circle*{2}}
      \put(50,20){\circle*{2}}
      \put(50,40){\circle*{2}}
      \put(-10,20){\line(0,1){20}}
      \put(-10,20){\vector(0,1){10}}
      \put(50,40){\line(0,-1){20}}
      \put(50,40){\vector(0,-1){10}}
      \qbezier(-10,40)(20,70)(50,40)
      \qbezier(-10,20)(20,-10)(50,20)
      \qbezier[40](-10,40)(20,30)(50,20)
      \qbezier[40](-10,20)(20,30)(50,40)
      \put(53,18){$c_j$}
      \put(53,38){$c_{j-1}$}
      \put(-25,18){$c_{i-1}$}
      \put(-19,38){$c_i$}
    \end{picture}
    \caption{}
    \label{fig:cyc}
  \end{center}
\end{figure}

If $(c_{j-1},c_{i-1})\in T$ then $(c_{i-1},c_i,\dots,c_{j-1})$ is an
ascending cycle in $T$ of length less than $m$, which contradicts
the requirement that $C$ is an ascending cycle on $T$ of minimal
length.  So $(c_{i-1},c_{j-1})\in T$. If $c_{i-1}<c_{j-1}$ then
$(c_{j-1},c_j,\dots,c_m,c_1,\dots,c_{i-1})$ is an ascending cycle in
$T$ of length less than $m$.  Hence, $c_{i-1}>c_{j-1}$.

Analogously, if $(c_i,c_j)\in T$ then
$(c_j,c_{j+1},\dots,c_p,c_1,\dots,c_i)$ is an ascending cycle in $T$
of length less than $m$.  So $(c_j,c_i)\in T$.  If $c_i>c_j$ then
$(c_i,c_{i+1},\dots,c_j)$ is an ascending cycle in $T$ of length less
than $m$. So $c_i<c_j$.

Now we have $c_{i-1}>c_{j-1}>c_j>c_i>c_{i-1}$, 
and we get an obvious contradiction.

We have shown that every minimal ascending cycle in $T$ is of
length $3$ or $4$ and thus have proved Theorem~\ref{th:semiacyclic}.  
\endproof

\subsection{The Orlik-Solomon algebra}
\label{sec:orlik-solomon}

In~\cite{OS} Orlik and Solomon gave the following combinatorial
description of the cohomology ring of an arbitrary hyperplane
arrangement.  Consider a complex arrangement $\cal A$ of affine
hyperplanes $H_1,H_2,\dots,H_N$ in the complex space $V\cong\C^n$ given by
\[
H_i:\ f_i(x)=0, \qquad i=1,\dots,N,
\]
where $f_i(x)$ are linear forms on $V$ (with a constant term).

We say that hyperplanes $H_{i_1},\dots,H_{i_p}$ are 
\emph{independent} if the codimension of the intersection
$H_{i_1}\cap\cdots\cap H_{i_p}$ is equal to $p$. 
Otherwise, the hyperplanes are \emph{dependent}.

Let $e_1,\dots,e_N$ be formal variables associated with the
hyperplanes $H_1,\dots,H_N$.  The \emph{Orlik-Solomon algebra}
$\mathrm{OS}(\A)$ of the arrangement $\A$ is generated over the
complex numbers by $e_1,\dots,e_N$ subject to the relations:

\begin{eqnarray}
\label{eq:OS1}
e_ie_j=-e_je_i,\qquad 1\leq i <j \leq N,
\\[.2in]
\label{eq:OS1.5}
e_{i_1}\cdots e_{i_p}=0,\qquad \textrm{if }H_{i_1}\cap\cdots\cap
H_{i_p}=\emptyset,\\[.2in]
\label{eq:OS2}
\sum_{j=1}^{p+1}(-1)^j\,e_{i_1}\cdots\widehat{e_{i_j}}\cdots
e_{i_{p+1}}=0,
\end{eqnarray}
whenever $H_{i_1},\dots,H_{i_{p+1}}$ are dependent.
(Here $\widehat{e_{i_j}}$ denotes that $e_{i_j}$ is missing.)

Let $C_\A=V-\bigcup_i H_i$ be the complement to the hyperplanes $H_i$ of $\A$,
and let $\mathrm{H}^*_{DR}(C_\A,\C)$ denote de Rham cohomology of
$C_\A$.
\begin{theorem}{\rm (Orlik, Solomon~\cite{OS})} \
\label{th:OS}
The map $\phi:\mathrm{OS}(\A)\to \mathrm{H}^*_{DR}(C_\A,\C)$
defined by
\[
\phi:\ e_i\mapsto\left[{df_i/f_i}\right]
\]
is an isomorphism.
\end{theorem}
Here $\left[{df_i/f_i}\right]$ is the cohomology class in
$\mathrm{H}_{DR}^*(C_\A,\C)$ of the differential form ${df_i/f_i}$.

We will apply Theorem~\ref{th:OS} to the Linial arrangement.  In this
case hyperplanes $x_i-x_j=1$, $i<j$, correspond to edges $(i,j)$ of
the complete graph~$K_n$. 

\begin{proposition}
  The Orlik-Solomon algebra $\mathrm{OS}(\L_{n-1})$ of the Linial
  arrangement is generated by $e_{vw}=e_{(v,w)}$, $1\le v<w\le n$
  subject to relations~(\ref{eq:OS1}),~(\ref{eq:OS1.5}),
  and also to  the following relations:
  \begin{equation}
    \label{eq:OS-linial}
    \begin{array}{c}
      e_{ab}e_{bc}e_{ac}-e_{ab}e_{bc}e_{cd}+e_{ab}e_{ac}e_{cd}-
      e_{bc}e_{ac}e_{cd}=0,\\[.1in]
      e_{ac}e_{bc}e_{bd}-e_{ac}e_{bc}e_{ad}+e_{ac}e_{bd}e_{ad}-
      e_{bc}e_{bd}e_{ad}=0.
    \end{array}
  \end{equation}
  where $1\leq a<b<c<d\leq n$ (cf.~Figure~\ref{fig:asc}).
\end{proposition}

\proof Let $C=(c_1,c_2,\dots,c_p)$ be a cycle in $K_n$. We say that
$C$ is \emph{balanced} if $\asc(C)=\des(C)$. We may assume that in
equation~(\ref{eq:OS2}) $i_1,i_2,\dots, i_p$ are edges of a balanced
cycle~$C$. We will prove~(\ref{eq:OS2}) by induction on~$p$. If $p=4$
then $C$ is of type~$C_1,C_2,C_3$, or~$C_4$
(see~Figure~\ref{fig:asc}). Thus $C$ produces one of the
relations~(\ref{eq:OS-linial}).  If $p>4$, then we can find $r\ne s$
such that both $C'=(c_r,c_{r+1},\dots,c_s)$ and
$C''=(c_s,c_{s+1},\dots,c_r)$ are balanced.  
Equation~(\ref{eq:OS2}) for $C$ is the sum of the equations for $C'$
and $C''$. Thus the statement follows by induction. \endproof

\begin{remark}
  {\rm This proposition is an analogue to the well-known description
    of the cohomology ring of the Coxeter arrangement~(\ref{eq:cox}),
    due to Arnold~\cite{Arnold}.  This cohomology ring is generated by
    $e_{vw}=e_{(v,w)}$, $1\leq v < w \leq n$, subject to
    relations~(\ref{eq:OS1}),~(\ref{eq:OS1.5}) and also the following 
    ``triangle'' equation:
    \[
    e_{ab}e_{bc}-e_{ab}e_{ac}+e_{bc}e_{ac}=0,
    \]
    where $1\leq a < b < c \leq n$.
    }
\end{remark}

\section{Truncated affine arrangements}
\label{sec:affine}
\neweq

In this section we study a general class of hyperplane arrangements
which contains, in particular, the Linial and Shi arrangements.

Let $a$ and $b$ be two integers such that $a+b\ge 2$.  Consider the
hyperplane arrangement $\A_{n-1}^{ab}$ in
$V_{n-1}=\{(x_1,\dots,x_n)\in\mathbb{R}^n\mid x_1+\cdots+x_n=0\}$
given by
\begin{equation}
\label{eq:affine_main}
x_i - x_j = -a+1,-a+2,\dots,b-1,\qquad 1\le i < j\le n.
\end{equation}

We call $\A_{n-1}^{ab}$ \emph{truncated affine arrangement} 
because it is a finite subarrangement of the affine arrangement 
of type $\widetilde{A}_{n-1}$ given by $x_i-x_j=k$, $k\in\mathbb{Z}$.

As we will see the arrangement $\A_{n-1}^{ab}$ has different behavior
in the \emph{balanced case} ($a=b$) and the \emph{unbalanced case} ($a\ne b$).

\subsection{Functional equations}

Let $f_n=f_{n}^{ab}$ be the number of regions of the arrangement
$\A_{n-1}^{ab}$, and let
\begin{equation}
  \label{eq:f}
  f(x)=\sum_{n\ge 0}f_n {x^n\over n!}
\end{equation}
be the exponential generating function for $f_n$.

\begin{theorem}
  \label{th:af.main} Suppose $a,b\ge 0$.
  \begin{enumerate}
  \item
    The generating function $f=f(x)$ satisfies the following
    functional equation:
    \begin{equation}
      \label{eq:funct_eq}
      f^{b-a}=e^{\textstyle x\cdot{f^{a}-f^{b}\over 1-f}}.
    \end{equation}
  \item
    If $a=b\geq 1$, then $f=f(x)$ satisfies the equation:
    \begin{equation}
      f=1+x\,f^{a},
      \label{eq:funct_eq:a=b}
    \end{equation}
  \end{enumerate} 
\end{theorem} 

Note that the equation~(\ref{eq:funct_eq:a=b}) can be formally
obtained from~(\ref{eq:funct_eq}) by l'H\^{o}pital's rule in the limit
$a\to b$.

In the case $a=b$ the functional equation~(\ref{eq:funct_eq:a=b})
allows us to calculate the numbers~$f_n^{aa}$ explicitly.  The following
statement was proved by P.~Headley~\cite{headley}.
\begin{corollary}
  \label{cor:a=b}
  The number $f_n^{aa}$ is equal to $an(an -1)\cdots(an-n +2)$.
\end{corollary}

The functional equation~(\ref{eq:funct_eq}) is especially simple in
the case $a=b-1$.  We call the arrangement~$\A_{n-1}^{a,a+1}$ the
\emph{extended Shi arrangement}.  In this case we get:

\begin{corollary}
\label{cor:gen-shi}
Let $a\geq 1$. The number $f_n$ of regions of the hyperplane
arrangement in $\mathbb{R}^n$ given by
\[
x_i-x_j=-a+1,-a+2,\dots,a,\qquad i<j,
\]
is equal to $f_n=(a\,n+1)^{n-1}$, and the exponential generating
function $f=\sum_{n\ge 0}f_n{x^n\over n!}$ satisfies the functional
equation $\displaystyle f=e^{x\cdot f^a}$.
\end{corollary}

In order to prove Theorem~\ref{th:af.main} we need several new
definitions.  A \emph{graded graph} is a graph $G$ on a set $V$ of
vertices labelled by natural numbers together with a function
$h:V\to\{0,1,2,\dots\}$, which is called a \emph{grading}. For $r\ge
0$ the vertices $v$ in $G$ such that $h(v)=r$ form the $r$th
\emph{level} of $G$.  Let $e=(u,v)$ be an edge in $G$, $u<v$.  We say
that the \emph{type} of the edge $e$ is the integer $t=h(v)-h(u)$ and
that a graded graph $G$ is \emph{of type} $(a,b)$ if the types of all
edges in $G$ are in the interval $[-a+1,b-1]=\{-a+1,-a+2,\dots,b-1\}$.

Choose a linear order on the set of all triples $(u,t,v)$, $u,v\in V$,
$t\in [-a+1,b-1]$.  Let $C$ be a graded cycle of type $(a,b)$. Every
edge $(u,v)$ in $C$ corresponds to a triple $(u,t,v)$, where $t$ is the
type of the edge $(u,v)$. Choose the edge $e$ in $C$ with the minimal
triple $(u,t,v)$.  We say that $C\setminus \{e\}$ is a \emph{broken
circuit of type} $(a,b)$.

Let $(F,h)$ be a graded forest.  We say that $(F,h)$ is
\emph{grounded} or that $h$ is a \emph{grounded} grading on the forest 
$F$ if each connected component in $F$ contains a vertex on the $0$th
level.

\begin{proposition}
\label{pr:af.broken_affine}
The number $f_n$ of regions of the arrangement~(\ref{eq:affine_main})
is equal to the number of grounded graded forests of type $(a,b)$ on
the vertices $1,2,\dots,n$ without broken circuits of
type~$(a,b)$.
\end{proposition}

\proof
By Corollary~\ref{cor:broken_graphic}, the number $f_n$ is equal to
the number of colored forests $F$ on the vertices $1,2,\dots,n$
without broken $A$-circuits.  Every edge $(u,v)$, $u<v$, in $F$ has a
color which is an integer from the interval $[-a+1,b-1]$. Consider the
grounded grading $h$ on $F$ such that for every edge $(u,v)$, $u<v$,
in $F$ of color~$t$ we have that $t=h(v)-h(u)$ is the type of $(u,v)$.  It
is clear that such a grading is uniquely defined. Then $(F,h)$ is a
grounded graded forest of type~$(a,b)$. Clearly, this gives a
correspondence between colored and graded forests. Then broken
$A$-circuits correspond to broken graded circuits. 
The proposition easily follows.
\endproof

From now on we fix the lexicographic order on triples $(u,t,v)$,
i.e., $(u,t,v)<(u',t',v')$ if and only if $u<u'$, or ($u=u'$ and
$t<t'$), or ($u=u'$ and $t=t'$ and $v<v'$). Note the order of $u$,
$t$, and $v$.  We will call a graded tree $T$ \emph{solid} if $T$
is of type $(a,b)$ and $T$ contains no broken circuits of type~$(a,b)$.

Let $T$ be a solid tree on $1,2,\dots,n$ such that 
vertex~$1$ is on the $r$th level.  If we delete the minimal vertex
$1$, then the tree $T$ decomposes into connected components
$T_1,T_2,\dots,T_m$.  Suppose that each component $T_i$ is connected
with $1$ by an edge $(1,v_i)$ where $v_i$ is on the $r_i$-th level.

\begin{lemma} 
\label{lem:adm} Let $T, T_1,\dots,T_m, v_1,\dots,v_m$, and
$r_1,\dots,r_m$ be as above.  The tree $T$ is solid if and only if
{\rm(a)} all $T_1,T_2,\dots,T_m$ are solid, {\rm(b)} for all i the
$r_i$-th level is the minimal nonempty level in $T_i$ such that
$-a+1\leq r_i-r\ge b-1$, and {\rm(c)} the vertex $v_i$ is the minimal vertex on
its level in $T_i$.  \end{lemma}

\proof First, we prove that if $T$ is solid then the conditions
(a)--(c) hold.  Condition~(a) is trivial, because if some $T_i$
contains a broken circuit of type~$(a,b)$ then $T$ also contains this
broken circuit.  Assume that  for some $i$ there is a vertex $v_i'$ on the
$r_i'$-th level in $T_i$ such that $r_i'<r_i$ and $r_i'-r\geq -a+1$.
Then the minimal chain in $T$ that connects  vertex $1$ with 
vertex $v_i'$ is a  broken circuit of type $(a,b)$.  Thus 
condition~(b) holds.  Now suppose that for some $i$  vertex $v_i$ is
not the minimal vertex $v_i''$ on its level.  Then the minimal chain
in~$T$ that connects  vertex~$1$ with~$v_i''$ is a broken circuit of
type~$(a,b)$.  Therefore,  condition~(c) holds too.

Now assume that conditions (a)--(c) are true.  We prove that $T$ is
solid.  For suppse not. Then $T$ contains a broken circuit
$B=C\setminus \{e\}$ of type~$(a,b)$, where $C$ is a graded circuit
and $e$ is its minimal edge.  If $B$ does not pass through vertex~$1$
then $B$ lies in~$T_i$ for some~$i$, which contradicts condition~(a).
We can assume that $B$ passes through vertex~$1$.  Since $e$ is the
minimal edge in $C$, $e=(1,v)$ for some vertex~$v'$ on level $r'$ in
$T$.  Suppose $v\in T_i$.  If $v'$ and $v_i$ are on different levels
in~$T_i$ then, by~(b), $r_i< r$.  Thus the minimal edge in~$C$
is~$(1,v_i)$ and not~$(1,v')$.  If $v'$ and $v_i$ are on the same
level in~$T_i$, then by~(c) we have $v_i< v'$.  Again, the minimal
edge in~$C$ is~$(1,v_i)$ and not~$(1,v')$.  Therefore, the tree~$T$
contains no broken circuit of type~$(a,b)$, i.e., $T$ is solid.
\endproof

Let $s_i$ be the minimal nonempty level in $T_i$, and let $l_i$ be the
maximal nonempty level in $T_i$. By Lemma~\ref{lem:adm}, the vertex~$1$
can be on the $r$th level, $r\in\{s_i-b+1,s_i-b+1,\dots,l_i+a-1\}$, and
for each such $r$ there is exactly one way to connect $1$ with $T_i$.



Let $p_{nkr}$ denote the number of solid trees (not necessarily
grounded) on the vertices $1,2,\dots,n$ which are located on
levels $0,1,\dots,k$ such that vertex~$1$ is on the $r$th level, $0\le
r\le k$.


Let 
\[ p_{kr}(x)=\sum_{n\ge 1}p_{nkr}{x^n\over n!},
\qquad p_k(x)=\sum_{r=0}^k p_{kr}(x).
\]


By the exponential formula (see~\cite[p.\ 166]{GJ}) and
Lemma~\ref{lem:adm}, we have 
\begin{equation}
p_{kr}'(x)=\exp(b_{kr}(x)),
\label{eq:2}
\end{equation}
where $b_{kr}(x)=\sum_{n\ge 1}b_{nkr}{x^n\over n!}$ and $b_{nkr}$ is
the number of solid trees $T$ on $n$ vertices located on the levels
$0,1,\dots, k$ such that at least one of the levels
$r-a+1,r-a+2,\dots,r+b-1$ is nonempty, $0\le r\le k$.  The polynomial
$b_{kr}(x)$ enumerates the solid trees on levels $1,2,\dots,k$
minus trees on levels $1,\dots,r-a$ and trees on levels $r+b,\dots,k$.
Thus we obtain
\[
b_{kr}(x)=p_k(x)-p_{r-a}(x)-p_{k-r-b}(x).
\]

By~(\ref{eq:2}), we get 
\[
p_{kr}'(x)=\exp(p_k(x)-p_{r-a}(x)-p_{k-r-b}(x)), 
\] 
where $p_{-1}(x)=p_{-2}(x)=\cdots=0$, $p_0(x)=x$, $p_k(0)=0$ for
$k\in\mathbb{Z}$. Hence
\[
p_k'(x)=\sum_{r=0}^k\exp(p_k(x)-p_{r-a}(x)-p_{k-r-b}(x)).
\]
Equivalently, 
\[
p_k'(x)\exp(-p_k(x))=\sum_{r=0}^k
\exp(-p_{r-a}(x))\, \exp(-p_{k-r-b}(x)).
\]
Let $q_k(x)=\exp(-p_k(x))$. We have
\begin{equation}
\label{eq:q_prime}
q_k'(x)=-\sum_{r=0}^k q_{r-a}(x)\,q_{k-r-b}(x),
\end{equation}
$q_{-1}=q_{-2}=\cdots=1$, $q_0=e^{-x}$, $q_k(0)=1$ for
$k\in\mathbb{Z}$.

The following lemma describes the relation between the polynomials
$q_k(x)$ and the number of regions of the arrangement ${\cal
A}_{n-1}^{ab}$.

\begin{lemma}
\label{lem:limit}
The quotient ${q_{k-1}(x)/ q_k(x)}$ tends to $\sum_{n\ge 0}f_n{x^n\over
n!}$ as $k\to\infty$.
\end{lemma}

\proof Clearly, $p_k(x)-p_{k-1}(x)$ is the exponential generating
function for the numbers of grounded solid trees of height less than
or equal to $k$. By the exponential formula (see~\cite[p.\ 166]{GJ})
$q_{k-1}(x)/ q_k(x)=\exp\left(p_k(x)-p_{k-1}(x)\right)$ is the
exponential generating function for the numbers of grounded solid
forests of height less than or equal to $k$.  The lemma obviously
follows from Proposition~\ref{pr:af.broken_affine}.  \endproof

All previous formulae and constructions are valid for arbitrary $a$
and $b$. Now we will take advantage of
the condition $a,b\ge 0$.  Let
\[
q(x,y)=\sum_{k\ge 0} q_k(x) y^k.
\]
By~(\ref{eq:q_prime}), we obtain the following differential equation
for $q(x,y)$:
\[
\begin{array}{rcl}
\displaystyle{\partial\over \partial x}\,
q(x,y)&=&-\left(a_y+y^{a}q(x,y)\right)\cdot
\left(b_y+y^{b}q(x,y)\right),\\[.2in]
q(0,y)&=&(1-y)^{-1},
\end{array}
\]
where $a_y:=(1-y^a)/(1-y)$.

This differential equation has the following solution:
\begin{equation}
\label{eq:qxt}
q(x,y)= {\textstyle b_y\,\exp(-x\cdot b_y)-a_y\,\exp(-x\cdot a_y)
\over\textstyle 
y^{a}\,\exp(-x\cdot a_y)-y^{b}\,\exp(-x\cdot b_y)}.
\end{equation}

Let us fix some small $x$.  Since $Q(y):=q(x,y)$ is an analytic
function of $y$, then $\gamma=\gamma(x)=\lim_{k\to\infty}{q_{k-1}/
  q_k}$ is the pole of $Q(y)$ closest to $0$ ($\gamma$ is the radius
of convergence of $Q(y)$ if $x$ is a small positive number).
By~(\ref{eq:qxt}), $\gamma^{a}\,\exp(-x\cdot
a_\gamma)-\gamma^{b}\,\exp(-x\cdot b_\gamma)=0$.  Thus, by
Lemma~\ref{lem:limit}, $f(x)=\sum_{n\ge 0}f_n{x^n\over n!}=\gamma(x)$
is the solution of the functional equation
\[
f^a\,e^{\textstyle -x\cdot{1-f^{a}\over 1-f}}= f^b\,e^{\textstyle
-x\cdot{1-f^{b}\over 1-f}},
\]
which is equivalent to~(\ref{eq:funct_eq}).

This completes the proof of Theorem~\ref{th:af.main}.
\endproof

\subsection{Formulae for the characteristic polynomial}
\label{sec:affine-formulae}

Let $\A=\A_{n-1}^{ab}$ be the truncated affine arrangement given
by~(\ref{eq:affine_main}).  Consider the characteristic polynomial
$\chi_n^{ab}(q)$ of the arrangement $\A_{n-1}^{ab}$.  Recall that
$\chi_n^{ab}(q)=q^{n-1}\poin_{\A_{n-1}^{ab}}(-q^{-1})$.

Let $\chi^{ab}(x,q)$ be the exponential generating function 
\[
\chi^{ab}(x,q)=1+\sum_{n>0}\chi_{n-1}^{ab}(q)\,\frac{x^n}{n!}\,.
\]
According to~\cite[Theorem~1.2]{S}, we have 
\begin{equation}
  \chi^{ab}(x,q)=f(-x)^{-q},
  \label{eq:exp}
\end{equation}
where $f(x)=\chi^{ab}(-x,-1)$ is the exponential generating
function~(\ref{eq:f}) for numbers of regions of~$\A_{n-1}^{ab}$.

Let $S$ be the \emph{shift operator} $S:f(q)\mapsto f(q-1)$. 
\begin{theorem}
\label{th:aff_poinc} 
Assume that $0\le a<b$.  Then
\[
\chi_n^{ab}(q)=(b-a)^{-n}(S^a+S^{a+1}+\cdots+S^{b-1})^n\cdot q^{n-1}.
\]
\end{theorem}

\proof The theorem can be deduced from Theorem~\ref{th:af.main}
and~(\ref{eq:exp}) (using, e.g., the Lagrange inversion formula).
\endproof

In the limit $b\to a$, using l'Hospital's rule, we obtain
\[
\chi^{aa}_n(q)=\left(S^a{\log S\over 1-S}\right)^n\cdot q^{n-1}.
\]

In fact, there is an explicit formula for~$\chi^{aa}(q)$. The
following statement easily follows from Corollary~\ref{cor:a=b} and
appears in \cite[??]{headley}\cite[proof of Prop.\ 3.1]{e-r}.
\begin{theorem}
  \label{th:q.a=b}
  We have 
  \[
  \chi_n^{a a}(q)=(q+1-an)(q+2-an)\cdots(q+n-1-an).
  \]
\end{theorem}

There are several equivalent ways to reformulate
Theorem~\ref{th:aff_poinc}, as follows:

\begin{corollary}
  \label{cor:aff-ways}
  Let $r=b-a$.  
\begin{enumerate}
\item We have
  \[
  \chi_n^{ab}(q)=r^{-n}\sum \left(q-\phi(1)-\cdots-\phi(n)\right)^{n-1},
  \]
  where the sum is over all functions $\phi:\,\{1,\dots,n\}\to\{a,\dots,b-1\}$.
\item We have
  \[
  \chi_n^{ab}(q)=r^{-n}\sum_{s,\,l\ge0} (-1)^l (q-s-an)^{n-1}{n\choose
    l}{s+n-rl-1\choose n-1}.
  \]
\item We have
  \[
  \chi_n^{ab}(q) = r^{-n}\sum {n\choose n_1,\dots,n_r} \left(q-an_1-\cdots
    -(b-1)n_r\right)^{n-1},
  \]
  where the sum is over all nonnegative integers $n_1,n_2,\dots,n_r$
  such that $n_1+n_2+\cdots+n_r=n$.
\end{enumerate}
\end{corollary}

\noindent
\textbf{Examples:}
\begin{enumerate}
\item ($a=1$ and $b=2$) \ The Shi arrangement ${\cal S}_{n-1}$ given
  by~(\ref{eq:1.shi}) is the arrangement~$\A_{n-1}^{12}$.  By
  Corollary~\ref{cor:aff-ways}.1, we get the following formula of
  Headley~\cite[??]{headley} 
  (generalizing the
  formula $r({\cal S}_{n-1})=(n+1)^{n-1}$ due to
  Shi~\cite[??]{shi}\cite{shi2}): 
  \[
  \chi_n^{1\,2}(q)=(q-n)^{n-1}.
  \]
\item ($a\geq 1$ and $b=a+1$) \ More generally, for the extended Shi
  arrangement ${\cal S}_{n-1,\,k}$ given by~(\ref{eq:1.ext.shi}), we
  have (cf.~Corollary~\ref{cor:gen-shi})
  \[
  \chi_n^{a,\,a{+}1}(q)=(q-an)^{n-1}.
  \]
\item ($a=0$ and $b=2$) \ In this case we get the Linial arrangement
  $\L_{n-1}={\cal A}_{n-1}^{02}$ (see Section~\ref{sec:linial}).  By
  Corollary~\ref{cor:aff-ways}.3, we have (cf.~Theorem~\ref{th:linial})
  \begin{equation}
    \label{eq:char_linial}
    \chi^{0\,2}_n(q) = 2^{-n}\sum_{k=0}^n {n\choose k}(q-k)^{n-1},
  \end{equation}
\item ($a\geq 0$ and $b=a+2$) \
  More generally, for the arrangement~$\A_{n-1}^{a,\,a+2}$, we have 
  \[
  \chi^{a,\,a+2}_n(q)= 2^{-n}\sum_{k=0}^n {n\choose k} (q-an-k)^{n-1}.
  \]
\end{enumerate}

Formula~(\ref{eq:char_linial}) for the characteristic polynomial
$\chi^{0\,2}_n(q)$ was earlier obtained by C.~Athanasiadis
\cite[Theorem~5.2]{Ath}.  He used a different approach based on a
combinatorial interpretation of the value of $\chi_n(q)$ for
sufficiently large primes~$q$.

[asymptotic behavior of $\chi_n^{ab}(q)$ --- to be inserted]

\subsection{Roots of the characteristic polynomial}
\label{sec:roots}

Theorem~\ref{th:aff_poinc} has one surprising  application
concerning the location of roots of the characteristic 
polynomial $\chi_n^{ab}(q)$

We start with the case $a=b$.  One can reformulate
Theorem~\ref{th:q.a=b} in the following way:

\begin{corollary}
  Let $a\geq 1$. The roots of the polynomial $\chi_n^{a a}(q)$ are the
  numbers $an-1,an-2,\dots,an-n+1$ (each with multiplicity $1$).  In
  particular, the roots are symmetric to each other with respect to
  the point $(2a-1)n/2$.
\end{corollary}

Now assume that $a\ne b $.  

\begin{theorem}
  \label{th:root-location}
  Let $a+b\geq 2$. All the roots of the characteristic polynomial
  $\chi_n^{ab}(q)$ of the truncated affine arrangement
  $\A_{n-1}^{ab}$, $a\ne b$, have real part equal to $(a+b-1)\,n/2$.
  They are symmetric to each other with respect to the point
  $(a+b-1)\,n/2$.
\end{theorem}

Thus in both cases the roots of the polynomial $\chi_n^{ab}(n)$
are symmetric to each other with respect to the point $(a+b-1)\,n/2$,
but in the case $a=b$ all roots are real, whereas in the case $a\ne b$
the roots are on the same vertical line in the complex plane
$\mathbb{C}$.  Note that in the case $a=b-1$ the polynomial
$\chi_n^{ab}(q)$ has only one root $an=(a+b-1)n/2$ with multiplicity
$n-1$.

The following lemma is implicit in a paper of Auric \cite{Auric} and
also follows from a problem posed by P\'olya \cite{Pol} and solved by
Obreschkoff \cite{Obr} (repeated in \cite[Problem V.196.1, pp.\ 70 and
251]{P-S}). For the sake of completeness we give a simple proof.

\begin{lemma} 
\label{lemma:realz}
Let $P(q)\in \mathbb{C}[q]$ have the property that every
root has real part $a$. Let $z$ be a complex number satisfying
$|z|=1$. Then every root of the polynomial $R(q) = (S+z)P(q) = P(q-1)
+ zP(q)$ has real part $a+\frac 12$.
\end{lemma}

\textbf{Proof.} We may assume that $P(q)$ is monic. Let
        $$ P(q) = \prod_j (q-a-b_j i),\ \ \ b_j\in\mathbb{R}. $$
If $R(w)=0$, then $|P(w)|=|P(w-1)|$. Suppose that $w=a+\frac 12+c+di$,
where $c,d\in\mathbb{R}$ and $i=\sqrt{-1}$. Thus
  $$ \left| \prod_j \left(\frac 12 + c +(d-b_j)i\right)\right| =
     \left|  \prod_j \left(-\frac 12 + c +(d-b_j)i\right)\right|. $$
If $c>0$ then $\left| \frac 12+c+(d-b_j)i\right| > \left| -\frac
12+c+(d-b_j)i\right|$. If $c<0$ then we have strict inequality in the
opposite direction. Hence $c=0$, so $w$ has real part $a+\frac
12$. \endproof 

\textbf{Proof of Theorem~\ref{th:root-location}.} All the roots of the
polynomial $q^{n-1}$ have real part 0. The operator
$T=(S^a+S^{a+1}+\cdots +S^{b-1})^n$ can be written as
  $$ T = S^{an}\prod_{j=1}^{b-1-a} (S-z_j)^n, $$
where each $z_j$ is a complex number of absolute value one (in fact, a
root of unity). The proof now follows from Theorem~\ref{th:aff_poinc}
and Lemma~\ref{lemma:realz}. \endproof

\subsection{Other root systems.}
\label{sec:conj}
The results of Subsections~\ref{sec:affine}.1--\ref{sec:affine}.3
extend, partly conjecturally, to all the other root
systems, as well as to the nonreduced root system $BC_n$ (the union of
$B_n$ and $C_n$, which satisfies all the root system axioms except the
axiom stating that if $\alpha$ and $\beta$ are roots satisfying
$\alpha=c\beta$, then $c=\pm 1$). Henceforth in this section when we
use the term ``root system,'' we also include the case $BC_n$.

Given a root system $R$ in $\mathbb{R}^n$ and integers $a$ and $b$
satisfying $a+b\geq 2$, we define the \emph{truncated $R$-affine
arrangement} ${\cal A}^{ab}(R)$ to be the collection of hyperplanes
  $$ \langle \alpha, x\rangle = -a+1, -a+2,\dots, b-1, $$
where $\alpha$ ranges over all positive roots of $R$ (with respect to
some fixed choice of simple roots). Here $\langle\ ,\ \rangle$ denotes
the usual scalar product on $\mathbb{R}^n$, and $x=(x_1,\dots,x_n)$.
As in the case $R=A_{n-1}$ we refer to the \emph{balanced case}
($a=b$) and \emph{unbalanced case} ($a\neq b$).

The characteristic polynomial for the balanced case was found by
Edelman and Reiner \cite[proof of Prop.\ 3.1]{e-r} for the root system
$A_{n-1}$ (see Theorem~\ref{th:q.a=b}), and conjectured (Conjecture
3.3) by them for other root 
systems. This conjecture was proved by Athanasiadis \cite[Cor.\ 7.2.3
and Thm.\ 7.7.6]{Athesis} for types $A,\ B,\ C,\ BC$, and $D$. For
types $A,\ B,\ C$ and $D$ the result is also stated in \cite[Thm.\
5.5]{Ath}.  We will not say anything more about the balanced case
here.

For the unbalanced case, we have considerable evidence (discussed
below) to support the following conjecture.

\begin{conjecture}
\label{conj:rpz}
Let $R$ be an irreducible root system in $\mathbb{R}^n$. Suppose that
the unbalanced truncated affine arrangement ${\cal A}={\cal
A}^{ab}(R)$ has $h({\cal A})$ hyperplanes. Then all the roots of the
characteristic polynomial $\chi_{\cal A}(q)$ have real part equal to
$h({\cal A})/n$.
\end{conjecture}

\textsc{Note.} (a) If all the roots of $\chi_{\cal A}(q)$ have the
same real part, then this real part must equal $h({\cal A})/n$, since
for any arrangement ${\cal A}$ in $\mathbb{R}^n$ the sum of the roots
of $\chi_{\cal A}(q)$ is equal to $h({\cal A})$.

(b) Conjecture~\ref{conj:rpz} implies the ``functional equation''
  \beq \chi_{\cal A}(q) = (-1)^n\chi_{\cal A}(-q+2h({\cal A})/n). 
  \label{eq:feq} \eeq
Thus $\chi_{\cal A}(q)$ is determined by around half of its
  coefficients (or values).

(c) Let $a+b\geq 2$ and $R=A_n,\ B_n,\  C_n$, or $D_n$. 
  Athanasiadis \cite[Thms.\ 7.2.1 and 7.2.4]{Athesis} has shown that
  except possibly when both $a=1$ and $R=C_n$ we have
  \beq \chi^{ab}_R(q) = \chi^{0,b-a}(q-ak), \label{eq:tranchi} \eeq
 where $k$ denotes the Coxeter number of $R$.
Presumably this equation also holds for the missing case $a=1$ and
  $R=C_n$. For $BC_n$ there is a similar result of Athanasidis
  \cite[Thm.\ 7.2.4]{Athesis}. 
These results and conjectures reduce Conjecture~\ref{conj:rpz} to the
case $a=0$ when $R$ is a classical root system. A similar reduction is
  likely to hold for the exceptional root systems.

(d) Conjecture~\ref{conj:rpz} is true for all the classical root
systems ($A_n,\ B_n,\ C_n,\ BC_n,\ D_n$). This follows from explicit
formulas found for $\chi^{ab}_R(q)$ by Athanasiadis \cite{ath3}
together with Lemma~\ref{lemma:realz}. The result of Athanasiadis is
the following.




\begin{theorem}
\label{thm:rootsys}

Up to a constant factor, we have the following characteristic
polynomials of the indicated arrangements. (If the formula has the
form $F(S)q^n$ or $F(S)(q-1)^n$, then the factor is $1/F(1)$.)
  $$ \begin{array}{rl} {\cal A}^{0,2k+2}(B_n): & (1+S^2+\cdots+
     S^{2k})^2 (1+S^2+\cdots+S^{4k+2})^{n-1} (q-1)^n\\[.1in]
     {\cal A}^{0,2k+2}(C_n): & \mbox{same as for ${\cal
        A}^{0,2k+2}(B_n)$}\\[.1in]
     {\cal A}^{0,2k+1}(B_n): & (1+S+\cdots+S^{2k})^2 (1+S^2+\cdots+
        S^{4k})^{n-1}q^n\\[.1in]
     {\cal A}^{0,2k+1}(C_n): & \mbox{same as for ${\cal
        A}^{0,2k+1}(B_n)$}\\[.1in]
     {\cal A}^{0,2k+2}(D_n): & (1+S^2)(1+S^2+\cdots+S^{2k})^4
        (1+S^2+\cdots+S^{4k+2})^{n-3}(q-1)^n \\[.1in]
     {\cal A}^{0,2k+1}(D_n): & (1+S+\cdots+S^{2k})^4 (1+S^2+\cdots+
        S^{4k})^{n-3}q^n\\[.1in]
     {\cal A}^{0,2k+2}(BC_n): & (1+S^2+\cdots+S^{2k})
        (1+S^2+\cdots+S^{4k+2})^n (q-1)^n\\[.1in]
     {\cal A}^{0,2k+1}(BC_n): &
        (1+S+\cdots+S^{2k})(1+S^2+\cdots+S^{4k})^n q^n.
   \end{array} $$
\end{theorem}

We also checked Conjecture~\ref{conj:rpz} for the arrangements ${\cal
A}^{02}(F_4)$ and ${\cal A}^{02}(E_6)$ (as well as the almost trivial
case ${\cal A}^{ab}(G_2),\ a\neq b$). The characteristic polynomials
are
  $$ \begin{array}{rl} {\cal A}^{02}(F_4): & q^4-24q^3+258q^2-1368q
     +2917\\[.1in]
   {\cal A}^{02}(E_6): & q^6-36q^5+630q^4-6480q^3+40185q^2-140076q
     +212002. 
  \end{array} $$ 

The formula for $\chi^{02}_{F_4}(q)$ has the remarkable alternative
form:
  $$ {\cal A}^{02}(F_4):\quad \frac 18 ((q-1)^4+3(q-5)^4+3(q-7)^4
     +(q-11)^4)-48. $$
Note that the numbers $1,5,7,11$ are the exponents of the root system
$F_4$. For $E_6$ the analogous formula is given by
  $$ {\cal A}^{02}(E_6): \quad \frac{1}{1008}P(q)-210, $$
where
  $$ P(q) = 61(q-1)^6+352(q-4)^6
   +91(q-5)^6+91(q-7)^6+352(q-8)^6+61(q-11)^6, $$
which is not as intriguing as the $F_4$ case.
It is not hard to see that the symmetry of the coefficient sequences
$(1,3,3,1)$ and $(61,352,91,91,352,61)$ is a consequence of equation
(\ref{eq:feq}) and the fact that if $e_1<e_2<\cdots<e_n$ are the
exponents of an irreducible root system $R$, then $e_i+e_{n+1-i}$ is
independent of $i$.

\end{document}